\documentclass[a4paper,12pt]{article}
\usepackage{amssymb,amsmath,amscd,amsfonts,amsthm,latexsym}
\usepackage{xcolor,graphpap,fancybox,framed,enumerate,graphicx}
\usepackage[english,activeacute]{babel}
\usepackage{needspace}
\usepackage[right=3cm,left=3cm,top=3cm,bottom=3cm]{geometry}
\usepackage{float}

\newtheorem{theorem}{Theorem}[section]
\newtheorem{prop}[theorem]{Proposition}
\newtheorem{coro}[theorem]{Corollary}
\newtheorem{lemma}[theorem]{Lemma}

\newenvironment{demo}{ \noindent \emph{\textbf{Proof:}}}{\hfill$\square$\\}

\newcommand{\RR}{\mathbb{R}}
\newcommand{\NN}{\mathbb{N}}
\newcommand{\CC}{\mathbb{C}}
\newcommand{\QQ}{\mathbb{Q}}
\newcommand{\TT}{\mathbb{T}}

\newcommand{\UU}{\mathbb{U}}

\newcommand{\Cc}{\mathcal{C}}

\newcommand{\Lc}{\mathcal{L}}

\newcommand{\Oc}{\mathcal{O}}

\newcommand{\grad}{\nabla}

\newcommand{\Un}{1\hspace{-1.5mm}1}

\newcommand{\dd}{\,{\text{\rm d}}}

\newcommand{\pc}{ \usefont{T1}{cmtl}{m}{n} \selectfont}

\newcommand{\ra}{\rangle}
\newcommand{\la}{\langle}
\newcommand{\G}{\Gamma}
\newcommand{\ddd}{\partial}
\newcommand{\iii}{{\, \vert\kern-0.25ex\vert\kern-0.25ex\vert\, }}
\newcommand{\ai}{{a_{\rm i}}}
\newcommand{\af}{{a_{\rm f}}}
\newcommand{\ui}{{u_{\rm i}}}
\newcommand{\uf}{{u_{\rm f}}}

\newdimen\texpscorrection
\newdimen\figcenter
\def\figurewithtex #1 #2 #3 #4 #5\cr{ Kll
  {\goodbreak\figcenter=\hsize\relax
  \advance\figcenter by -#4truecm
  \divide\figcenter by 2
  \begin{figure}[hbt]
  \vskip #3truecm\noindent\hskip\figcenter
  \includegraphics{#1}{\hskip\texpscorrection\input #2 }
  \vskip 0.8truecm{\baselineskip=0.8\baselineskip
  \noindent \vbox{\noindent {\footnotesize #5}}\par}
  \end{figure}}}
\def\point#1 #2 #3 {\rlap{\kern #1 truecm
\raise #2 truecm \hbox{#3}}}

\numberwithin{equation}{section}

\begin{document}

\title{\bf Permuting quantum eigenmodes by a quasi-adiabatic motion of a potential wall}

\author{Alessandro \textsc{Duca}\footnote{Universit\'e 
Grenoble Alpes, CNRS, Institut Fourier, F-38000 Grenoble, France email: {\pc 
alessandro.duca@univ-grenoble-alpes.fr}} ,  Romain 
\textsc{Joly}\footnote{Universit\'e Grenoble Alpes, CNRS, 
Institut Fourier, F-38000 Grenoble, France email: {\pc 
romain.joly@univ-grenoble-alpes.fr}} {~\&~}
Dmitry \textsc{Turaev}
\footnote{Imperial College, London SW7 2AZ, UK, and Higher School of Economics - Nizhny
Novgorod Lobachevsky University and Higher School of Economics of Nizhny Novgorod, Bolshaya Pecherskaya 25/12, 603155, Nizhny Novgorod, Russia, email: {\pc d.turaev@imperial.ac.uk }}}
\date{}

\maketitle

\begin{abstract} We study the Schr\"odinger equation $i\ddd_t\psi=-\Delta\psi+V\psi$ on $L^2((0,1),\CC)$ where $V$ is a very high and localized potential wall. We aim to perform permutations of the eigenmodes and to control the solution of the equation. We consider the process where the position and the height of the potential wall change as follows. First, the potential increases from zero to a very large value, so a narrow potential wall is formed that almost splits the interval into two parts; then the wall moves to a different position, after which the height of the wall decays to zero again. We show that even though the rate of the variation of the potential's parameters can be arbitrarily slow, this 
process alternates adiabatic and non-adiabatic dynamics, leading to a non-trivial permutation of the eigenstates. Furthermore, we consider potentials with several narrow walls and we show how an arbitrarily slow motion of the walls can lead the system from any given state to an arbitrarily small neighborhood of any other state, thus proving the approximate controllability of the above Schr\"odinger equation by means of a soft, quasi-adiabatic variation of the potential.

\vspace{3mm}

\end{abstract}


\section{Introduction}\label{intro}

In this paper, we control the eigenstates of the Schr\"odinger equation in a bounded interval by moving very slowly a sharp and narrow potential wall. Our aim is twofold. First, we exhibit how to \emph{permute eigenstates through slow cyclic processes} which seem to violate the adiabatic principle. Second, we provide a \emph{new method for the approximate control} of the Schr\"odinger equation.    

\vspace{3mm} 

One of the basic principles of quantum mechanics is that quantum numbers are preserved when
parameters of the system change sufficiently slowly \cite{BF,Kato,AE}. That is, a slow variation of the Hamiltonian operator of a quantum-mechanical system makes it evolve adiabatically: if one prepares the initial state with a definite energy and starts to change the Hamiltonian slowly, the system will stay close to the
state of the specific energy, defined by the same set of quantum numbers, for a very long time. In the absence of symmetries, this means that if we order the eigenstates of the Hamiltonian by the increase
of the energy, then the adiabatic evolution along a generic path from a Hamiltonian $H_1$ to
a Hamiltonian $H_2$ must lead the system from the $k$-th eigenstate of $H_1$ to the $k$-th eigenstate of $H_2$, with the same $k$. In particular, a slow cyclic (time-periodic) variation of the Hamiltonian is expected
to return the system to the initial eigenstate (up to a phase change) after each period, for many periods.

In this paper, we describe a class of (quasi)adiabatic processes which violate this principle in the following sense. We consider a particle in a bounded region, subject to an external potential. We show how an arbitrarily slow cyclic change in the potential can lead the system close to a different eigenstate after each  cycle. Moreover, the evolution over the cycle of our process may act as any finite permutation of the energy eigenstates.  By building further on this idea, we are able to perform a global approximate control of the Schr\"odinger equation using adiabatic arguments: by a slow change of potential we can lead the system from any given state to any other one, up to an arbitrarily small error.

The results provide a rigorous implementation of a general idea from \cite{Dmitry}. Consider a periodic
and slow variation of a Hamiltonian such that for a part of the period the system acquires an additional
quantum number (a quantum integral) that gets destroyed for the rest of the period. Then the system will
evolve adiabatically, however it can find itself close to a different eigenstate at the end of each period.
In this way, the periodic creation and destruction of an additional quantum integral gives rise to much more general and 
diverse types of adiabatic processes than it was previously thought.

\vspace{3mm}

Here, we consider the one-dimensional Schr\"odinger equation in $L^2((0,1),\CC)$ with a locally supported potential with time-dependent position and height: 
\begin{equation}\label{eq}\tag{SE}
\begin{cases} i\ddd_t u(t,x)=-\ddd_{xx}^2 u(t,x)+ V(t,x)u(t,x),\ \ \ \ \ \ \ & x\in (0,1),~t > 0 \\
u(t,0)=u(t,1)=0, & t\geq 0,\\
u(t=0,\cdot)=u_0(\cdot) \in L^2((0,1),\CC).\end{cases}
\end{equation}
We define the potential $V$ as follows. Let $\rho\in\Cc^2(\RR,\RR^+)$ be a 
non-negative function with support $[-1,1]$ such that $\int_\RR \rho(s) \dd s=1$. We consider $\eta\in\Cc^2(\RR,\RR^+)$, $I\in\Cc^2(\RR,\RR^+)$ and $a\in\Cc^2(\RR,(0,1))$.
We set
\begin{equation}\label{eq_potentiel}
V(t,x)~=~ I(t)\,\rho^{\eta(t)}(x-a(t))~~~~~\text{ where }~~~\rho^\eta=\eta\,\rho(\eta\,\cdot\,)~. 
\end{equation}
We notice that \eqref{eq} is a free Schr\"odinger equation when $\eta\equiv 0$ or $I\equiv 0$, while $\rho^\eta$ is close to the Dirac delta-function for large $\eta$. When $I$ and $\eta$ are very large, the potential $V=I\rho^\eta$ describes a very high and thin wall that splits
the interval $(0,1)$ into two parts $(0,a)$ and $(a,1)$. The parameter $a$ enables to move the location of the wall, while $I$ controls 
its height and $\eta$ defines its sharpness. We are interested in controlling the eigenmodes of the Schr\"odinger equation by suitable motions of $\eta$, $I$ and $a$.

\vspace{3mm}

The article is inspired by the work \cite{Dmitry} where the Schr\"odinger equation was considered in domains which are periodically divided into disconnected parts. In particular, one can consider the following example:
\begin{equation}\label{eq_dmitry}
\begin{cases} i\ddd_t u(t)=-\ddd_{xx}^2 u(t),&  x\in (0,a)\cup(a,1),\\
u(t,0)=u(t,1)=0,  & t\geq 0,\\
u(t,a-0)=u(t,a+0),& t\geq 0,\\
\alpha u(t,a)+(1-\alpha)(\ddd_x u(t,a+0)-\ddd_x u(t,a-0))=0 & t\geq 0,\\
u(t=0,\cdot)=u_0(\cdot) \in L^2((0,a)\cup(a,1),\CC).\end{cases}
\end{equation}
Such equation corresponds to \eqref{eq} when $V$ is a singular
Dirac potential $I \delta_{x=0}$ with $I=\alpha/(1-\alpha)$, and $\alpha\in [0,1]$ is slowly changing, continuous function of time. It is assumed that $\alpha\equiv 0$ at $t$ close to $0$ and $T$ and $\alpha\equiv 1$ for some closed subinterval of $(0,T)$. By construction, this equation aims to describe adiabatic transitions from a free particle confined in the interval $(0,1)$ (when $\alpha=0$) to a free particle in the pair of disjoint intervals $(0,a)$ and $(a,1)$ (when $\alpha=1$).  In \cite{Dmitry}, it is shown that the adiabatic evolution starting at $t=0$ with an eigenstate of the Laplacian with the Dirichlet boundary conditions on the interval $(0,1)$ leads this system, typically, to the vicinity of a different
eigenstate at $t=T$, and that repeating the same process for many time periods makes the energy grow exponentially in time. In the present
paper we do not discuss the phenomenon of the exponential energy growth. Instead, we analyze in detail the evolution over one period and, in particular, pursue a new idea of controlling the evolution of the system by controlling the slow rate with which parameters of the system are varied.

The analysis of \cite{Dmitry} is not rigorous because of the singularity of the Dirac potential. We therefore replace the singular model (\ref{eq_dmitry}) by the well-posed time-dependent PDE \eqref{eq}. Such choice helps us to characterize the dynamics of the equation rigorously, even though it lets the tunneling effect appear since the potential is not a perfect infinite wall. Despite of this tunneling effect, we can show that the evolution can mimic the permutations of eigenmodes described in \cite{Dmitry}. Moreover, the tunneling effect can also be helpful to obtain the global control of the equation by varying the speed with which the potential wall moves. 

\vspace{5mm}

\begin{figure}[ht]
\begin{center}
\resizebox{0.5\textwidth}{100pt}{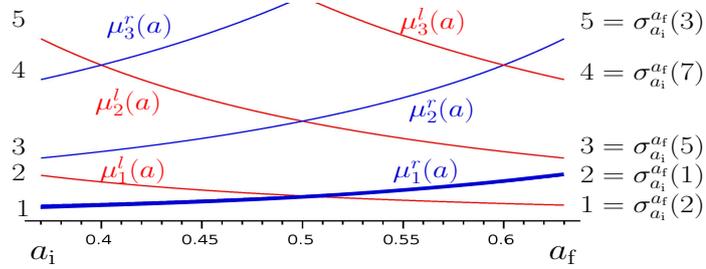}
\end{center}
\caption{\it The figure represents an example of quasi-adiabatic permutation $\sigma_\ai^\af$. The eigenvalues of the Dirichlet Laplacian on $(0,a)\cup(a,1)$ with $a\in[a_i,a_f]$ are decomposed into two sequences given by \eqref{def_mu}. When $a$ goes from $\ai$ to $\af$, the eigenvalues corresponding to the different sequences can cross, which yields $\sigma_\ai^\af$. The line of eigenvalues $\mu_1^r(a)$ corresponds to the evolution of the energy value when the position of the potential wall moves from  $a=\ai$ to $a=\af$,
as shown in the middle of Figure \ref{fig_th}, enabling to transform the first eigenmode into the second one.}\label{fig-expli-permut}
\end{figure}

{\noindent \bf \underline{Main results: Permutation of the ei}g\underline{enstates}}\\[1mm] 
\noindent We start by formally introducing the permutation of eigenmodes that was described in \cite{Dmitry}. 
Let $a\in (0,1)$ and let $\NN^*$ be
the set of strictly positive integers. We set
\begin{equation}\label{def_mu}
\mu_p^l(a)=\frac{p^2\pi^2}{a^2}~~~\text{ and }~~~\mu_q^r(a)=\frac{q^2\pi^2}{(1-a)^2}~~~\text{ with }~~~p,q\in\NN^*.
\end{equation}
These numbers are the eigenvalues of the Dirichlet-Laplacian operator in the left and the right parts of the split interval $(0,a)\cup(a,1)$. For any initial and final positions $\ai$ and $\af$ in $(0,1)\setminus\QQ$, we define the {\em quasi-adiabatic permutation} $\sigma_\ai^\af:\NN^*\rightarrow\NN^*$ as follows.
Since $\ai$ is irrational, it follows that $\mu_p^l\neq \mu_q^r$ for all $p,q\in\NN^*$, so we can order
$$\{\mu_p^l(\ai)\}_{p\in\NN^*}\cup\{\mu_q^r(\ai)\}_{q\in\NN^*}$$ as a strictly increasing sequence. Let $k\in\NN^*$. The $k-$th element of this sequence is either a number $\mu_{p_0}^l(\ai)$ or a number $\mu_{q_0}^r(\ai)$ for some $p_0$ or $q_0$. When $a$ changes from $\ai$ to $\af$, the corresponding eigenvalue $\mu_{p_0}^l(\ai)$ or $\mu_{q_0}^r(\ai)$ smoothly goes to $\mu_{p_0}^l(\af)$ or, respectively, to $\mu_{q_0}^r(\af)$, with the same $p_0$ (or $q_0$). As $\af\not\in\QQ$, the union
$$\{\mu_p^l(\af)\}_{p\in\NN^*}\,\cup\,\{\mu_q^r(\af)\}_{q\in\NN^*}$$
can, once again, be ordered as a strictly increasing sequence in a unique way. Thus, we define $\sigma_\ai^\af(k)$ as the number such that 
$\mu_{p_0}^l(\af)$ (or $\mu_{q_0}^r(\af)$) is the $\sigma_\ai^\af(k)-$th element of the sequence obtained by the ordering of $\{\mu_p^l(\af)\}_{p\in\NN^*}\,\cup\,\{\mu_q^r(\af)\}_{q\in\NN^*}$. Figure \ref{fig-expli-permut} illustrates an example of this permutation.

\vspace{3mm}

The first main result of this paper is given by the following theorem. 
\begin{theorem}\label{th}
Let $\ai,\af\in(0,1)\setminus\QQ$, $N\in\NN^*$, $\varepsilon>0$ and $\kappa>0$. There exist $T>0$,
\begin{itemize}
\item $\eta\in\Cc^\infty([0,T],\RR^+)$ with $\|\eta'\|_{L^\infty([0,T],\RR)}\leq \kappa$,
\item $I\in\Cc^\infty([0,T],\RR^+)$ with $\|I'\|_{L^\infty([0,T],\RR)}\leq \kappa$ and $I(0)=I(T)=0$,
\item $a\in\Cc^\infty ([0,T],(0,1))$ with $\|a'\|_{L^\infty([0,T],\RR)}\leq \kappa$, 
$a(0)=\ai$ and $a(T)=\af,$ 
\end{itemize}
such that the evolution defined by the linear Schr\"odinger equation \eqref{eq} with the potential $V$ given by \eqref{eq_potentiel} realizes the quasi-adiabatic permutation $\sigma_\ai^\af$. Namely, let $\Gamma_{s}^{t}$ be the unitary propagator generated in the time interval $[s,t]\subset[0,T]$ by equation \eqref{eq} with the potential \eqref{eq_potentiel}. Then, for all $k\leq N$, there exist $\alpha_k\in\CC$ with 
$|\alpha_k|=1$ such that 
$$\big\|~\Gamma_{0}^{T}\sin(k\pi x)\,-\,\alpha_k\sin(\sigma_\ai^\af(k)\pi x)~\big\|_{L^2}\leq \varepsilon~. $$
\end{theorem}
\begin{figure}[ht]
\begin{center}
\resizebox{\textwidth}{!}{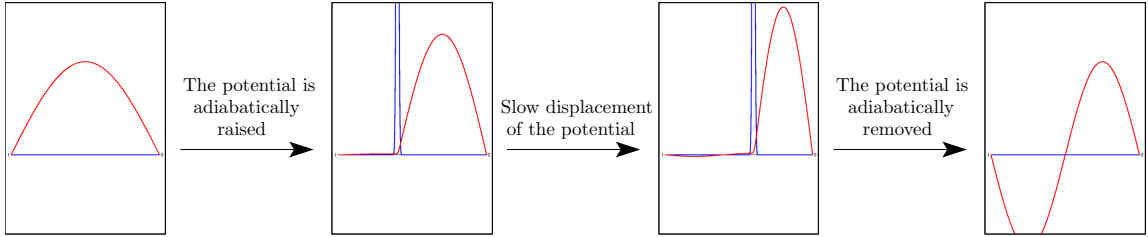}
\end{center}
\caption{\it The figure represents a control path transforming the first mode to the second one, according to Theorem \ref{th}. The change in energy when the potential wall moves to the right follows the line of eigenvalues $\mu_1^r(a)$ shown in Figure \ref{fig-expli-permut}.
}\label{fig_th}
\end{figure}
The trajectory of the potential wall $V$ given by Theorem \ref{th} is as follows (see Figure \ref{fig_th}). First, we slowly grow a very thin  
potential wall at $a=\ai$, until its height $I$ reaches a sufficiently large value so, in a sense, the interval gets almost split into two parts. During this stage the system evolves adiabatically: starting with the $k$-th eigenstate $\sin(k\pi x)$ of the Dirichlet Laplacian on $(0,1)$ it arrives close to the $k$-th eigenstate of the Dirichlet Laplacian on $(0,\ai)\cup(\ai,1)$. The latter eigenstates are localized either in $(0,\ai)$ or in $(\ai,1)$ (see Figure \ref{fig_th}), so the particle gets almost completely localized in one of these intervals (if the wall's height $I$ gets sufficiently high). After that, we slowly change the wall's position $a$. When the $k-$th eigenvalue of the Dirichlet Laplacian on $(0,a)\cup(a,1)$ is a double eigenvalue, a fully adiabatic evolution would lead to a strong tunneling so we must leave the adiabatic strategy to ensure that the particle remains trapped in $(0,a)$ or $(a,1)$ and that the system stays close to the corresponding ``left'' eigenstate (with the eigenvalue $\mu_{p_0}^l(a)$) or ``right'' eigenstate (with the eigenvalue $\mu_{q_0}^r(a)$). Thus, at the end of this stage the system will arrive close to the $\sigma_\ai^\af(k)-$th eigenstate of the Dirichlet Laplacian on $(0,\af)\cup(\af,1)$. Finally, we adiabatically decrease the potential $V$ to zero - and the system finds itself close to the eigenstate $\sin(\sigma_\ai^\af(k)\pi x)$. 

The main mathematical difficulty of this article occurs when the tunneling effect becomes strong and when we must thus leave the adiabatic regime. To avoid the tunneling, we need to accelerate but a fast motion of the high and sharp potential wall potentially generates large unstability. The main trick consists in controlling the error terms due to this short non-adiabatic motion.

\vspace{3mm}

\begin{figure}[ht]
\centering
\includegraphics[width=0.5\textwidth,height=100pt]{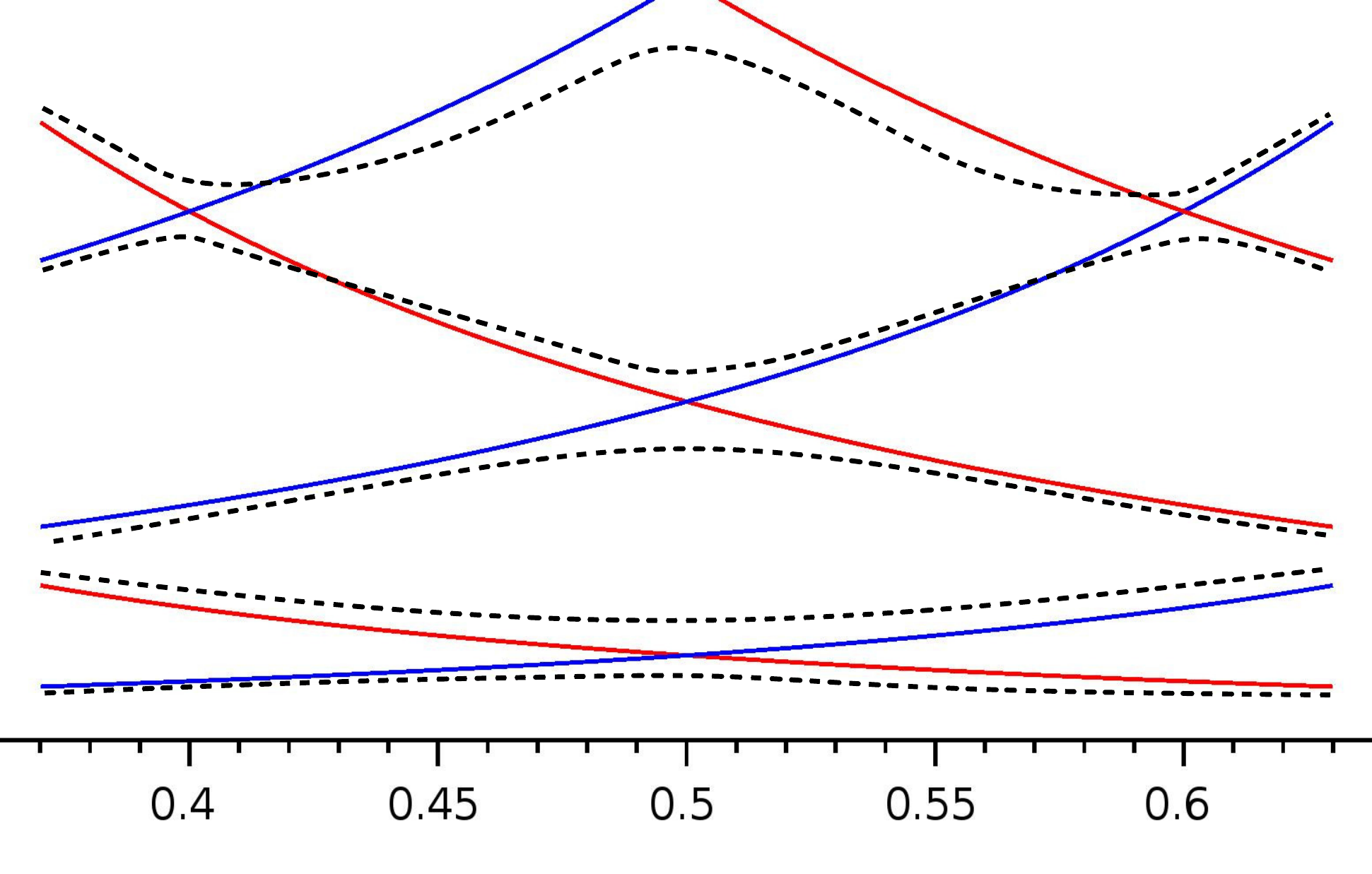}
\caption{\it The dashed lines represent the eigenvalues of the Hamiltonian $-\ddd_{xx}^2+V(t,\cdot)$ as functions of $t\in [0,T]$ with the potential given by given by \eqref{eq_potentiel} with sufficiently large $I$ and $\eta$. These curves are close to the set
of lines of eigenvalues of the singular problem \ref{eq_dmitry} shown in Figure \ref{fig-expli-permut}, but they never intersect. 
As a consequence, a fully adiabatic motion would preserve the ordering of the eigenmodes in this framework.
 }\label{approximating_curves}
\end{figure}

We stress that the dependence with respect to the time in the potential $V(t,x)$ constructed in Theorem \ref{th} can be as slow as desired (the constant $\kappa$ that bounds the speed of the parameter change can be as small as we need). This means the control we apply to the system in order to achieve the
permutation of the eigenstates is soft. However, it would be wrong to think of this process as fully adiabatic (we use the term ``quasi-adiabatic'' instead). Indeed, fully adiabatic dynamics would preserve the ordering of the eigenmodes since the spectrum of the Hamiltonian $-\ddd_{xx}^2+V(t,\cdot)$ given by \eqref{eq_potentiel} is simple for every $t\in [0,T]$, as presented in Figure \ref{approximating_curves} (see Proposition \ref{prop-simple}). Would the speed of parameters' change be too slow, the evolution would trace the eigenstates of the instantaneous Hamiltonian so tightly that no permutation of the eigenstates would be possible. Contrary to that, the permutation of eigenstates is achieved in Theorem \ref{th} by choosing the speed of the potential wall to be much faster than the rate of tunneling, that occurs when the potential wall approaches the values of $a$ where different eigenvalues become close. These moments correspond to the crossings of 
left and right eigenvalues $\mu_p^l(a)$ and $\mu_q^r(a)$ of the singular problem \eqref{eq_dmitry} (e.g. $a=1/2$ in Figure \ref{fig_th}). 
As we see, due to the presence of a small gap between the corresponding eigenvalues of the non-singular problem 
\eqref{eq}, the quasi-adiabatic evolution at these moments is controlled by the rate of the parameters' variation. This key fact is used in the next theorem, which provides a more general result than Theorem \ref{th}.

\vspace{3mm}

Note that there is enough freedom in this construction. In fact, the class of controls that can be used in order to achieve the claim of Theorem \ref{th} is quite wide, as one can see in the proof (see the discussion in Section \ref{applications}). The robustness of the proposed control is an important aspect of our results. In particular, $\eta(t)$ can be taken constant; we can also make $I(t)$ identically zero on some time intervals around the end points $t=0$ and $t=T$ of the control interval. The later means the process can be repeated periodically, e.g. realizing the exponential heating process described in \cite{Dmitry} (see Section \ref{energy_growth}). We also remark that, even though the phase shifts appearing in Theorem \ref{th} are not relevant from a physical point of view, they can be easily removed in order to obtain $\alpha_k=1$ for all $k$ (we postpone further explanations to Section \ref{section_phase} since its proof differs from the spirit of our core strategy).

\vspace{5mm}

{\noindent \bf \underline{Main results: Control of the Schr\"odin}g\underline{er e}q\underline{uation}}\\[1mm] 
\noindent Let us consider a more general class of potentials $V(t,x)$ with several slowly moving walls. Namely, let
\begin{equation}\label{eq_potentiel_2}
V(t,x)~=~\sum_{j=1}^J  I_j(t) \rho^{\eta_j(t)}(x-a_j(t))
\end{equation}
with $\{\eta_j\}_{j\leq J},\{I_j\}_{j\leq J}\subset\Cc^\infty([0,T],\RR^+)$ and $\{a_j\}_{j\leq J}\subset\Cc^\infty ([0,T],(0,1))$. Taking the 
number $J$ of the walls large enough and tuning the rate with which the parameters of the potential vary with time allows us to obtain the following {\em  approximate controllability} result.
\begin{theorem}\label{th3} 
Let $\varepsilon>0$, $\kappa>0$ and let $u_{\rm i}$ and $u_{\rm f}$ be functions from $L^2((0,1),\CC)$ with $\|u_{\rm i}\|_{L^2}=\|u_{\rm f}\|_{L^2}$. There exist $J\in\NN$, $T>0$, and smooth functions $\{\eta_j\}_{j\leq J},$  $\{I_j\}_{j\leq J}\subset\Cc^\infty([0,T],\RR^+)$, and $\{a_j\}_{j\leq J}\subset\Cc^\infty ([0,T],(0,1))$, with time derivatives bounded by $\kappa$ in the absolute value, such that
$$\|\Gamma_0^T u_{\rm i}\,-\,u_{\rm f}\|_{L^2}~\leq~\varepsilon~,$$
where $\Gamma_{s}^{t}$ denotes the unitary propagator generated by the linear Schr\"odinger equation \eqref{eq} in the time interval $[s,t]\subset[0,T]$ with the potential $V$ given by \eqref{eq_potentiel_2}.
\end{theorem}

Like in theorem Theorem \ref{th}, the functions $\eta_j$ can be taken constant, and the functions $I_j$ can be taken identically zero 
at the beginning and the end of the control interval $[0,T]$ (so the process described by Theorem \ref{th3} can be repeated several times). 
Theorem \ref{th3} is stronger than Theorem \ref{th}, however we formulated Theorem \ref{th} separately, as it provides a simpler
control protocol. In particular, the necessary number $J$ of walls in potential \eqref{eq_potentiel_2} depends on $u_{\rm i}$ and $u_{\rm f}$
and on the accuracy $\varepsilon$.

\vspace{5mm}

{\noindent \bf \underline{Some biblio}g\underline{ra}p\underline{h}y}\\[1mm] 
\noindent A peculiarity of our work is the nature of the control strategy. Adiabatic controls are not so common in literature (see for instance \cite{ugo2,CT04,CT06}), controls that alternate adiabatic and non-adiabatic motions are even more unusual. In this sense, even though the approximate controllability was already known for \eqref{eq}, our method provides a new and different way to perform it.

\vspace{1mm}

The controllability of Schr\"odinger equations of the form \eqref{eq} is commonly studied with potentials of the form $V(t,\cdot)=v(t) B$ with $t\in[0,T]\subset\RR^+$. The function $v$ represents the time-dependent intensity of a controlling external field described by the bounded symmetric operator $B$. Equation of such type, called bilinear Schr\"odinger equation, are known to be not exactly controllable in $L^2((0,1),\CC)$ when $v\in L^{r}_{loc}(\RR^+,\RR)$ with $r>1$, see the work \cite{ball} by Ball, Mardsen, and Slemrod.  The turning point for this kind of studies is the idea of controlling the equation in suitable subspaces of $L^2((0,1),\CC)$ introduced by Beauchard in \cite{be1}. Following this approach, several works achieved exact controllability results for the bilinear Schr\"odinger equation, see e.g. \cite{laurent,mio1,mio2,morgane1,morganerse2}. The global approximate controllability of bilinear quantum systems has been proven with the help of various techniques. We refer to \cite{milo,nerse2} for Lyapunov techniques, while we cite \cite{ugo2,ugo3} for adiabatic arguments, and \cite{nabile,ugo} for Lie-Galerkin methods. 

\vspace{1mm}

The controllability of linear Schr\"odinger equations with controls on the boundaries has been established e.g. in \cite{lion,cagn} by the multiplier method, in \cite{rauch,burqa,lollo} by
microlocal analysis, and in \cite{44,kaska,330} with the use of Carleman estimates. Concerning the controllability of various PDEs by moving boundaries, we refer to the works \cite{kais2,kais1,kar2,kar1,duclos}.

\vspace{5mm}

{\noindent \bf \underline{Scheme of the article}}\\[1mm]
\noindent In Section \ref{preliminaries}, we discuss some basic properties of the Schr\"odinger equation \eqref{eq}. In particular, Section \ref{well-posedness} ensures the well-posedness of the equation, while Section \ref{section_adiab} establishes the validity of a classical adiabatic result. 

\noindent
In Section \ref{section_spectre}, we study the spectrum of the Hamiltonian in equation \eqref{eq}. We show that, for $\eta$ and $I$ sufficiently large, such Hamiltonian approximates, from the spectral point of view, the Dirichlet Laplacian in the split segment considered in \cite{Dmitry}.

\noindent
In Section \ref{section_fast_translation}, we prove that the evolution defined by \eqref{eq} does not change too much the 
states which are small in the support of $V$ when we accelerate the motion of the potential wall for a short time.

\noindent
In Section \ref{section-proof-of-theorem}, we provide the proof of Theorem \ref{th} by gathering the results from the previous sections.

\noindent
In Section \ref{further_discussions} discuss further applications of the techniques behind Theorem \ref{th}. 
In particular, we provide the proof of Theorem \ref{th3}.

\vspace{3mm}

{\noindent \bf Acknowledgements:} the present work has been initiated by a discussion in the fruitful atmosphere of the Hale conference in Higashihiroshima. This work was supported by RSF grant 19-11-00280, RSF grant 19-71-10048 and the project {\it ISDEEC} ANR-16-CE40-0013. The third author also acknowledges support by EPSRC and Russian Ministry of Science and Education project 2019-220-07-4321.


\section{Preliminaries}\label{preliminaries}

\subsection{Basic notations and results}\label{well-posedness}
We equip the Hilbert space $L^2((0,1),\CC)$ with the scalar product 
$$\la\,\psi_1\,|\,\psi_2\,\ra_{L^2}=\int_0^1\overline{\psi_1(x)}\psi_2(x)\dd x,\ \ \ \ \ \forall\psi_1,\psi_2\in L^2((0,1),\CC)$$
and the corresponding norm $\|\cdot\|_{L^2}=\sqrt{\la\,\cdot\,|\,\cdot\,\ra_{L^2}}$. We consider the classical Sobolev's spaces $H^m((0,1),\CC)$ for $m\geq 0$ with the standard norms, the space $H^1_0((0,1),\CC)=\{u\in H^1((0,1),\CC),~u(0)=u(1)=0\}$ and the space $H^{-1}((0,1),\CC)$, which is the dual space of $H^1_0((0,1),\CC)$ with respect to the $L^2-$duality. 

\vspace{3mm}

Defining solutions of an evolution equation with a time-dependent family of operators is nowadays a classical result (see \cite{Tanabe}). In our case, we deduce the well-posedness of equation \eqref{eq} from the results in \cite{Kisynski} (see also \cite{Teufel}).
\begin{theorem}\label{th-Kisynski}
{\bf (Kisy\'nski, 1963)}\\
Let $X$ be a Hilbert space and let $(H(t))_{t\in [0,T]}$ be a family of self-adjoint positive operators 
on $X$ such that $X^{1/2}=D(H(t)^{1/2})$ is independent of time $t$. Also set 
$X^{-1/2}=D(H(t)^{-1/2})=(X^{1/2})^*$ and assume that $H(t):X^{1/2}\rightarrow 
X^{-1/2}$ is of class $\Cc^2$ with respect to $t\in [0,T]$. Assume that
there exists $\gamma>0$ and $\kappa\in\RR$ such that,
\begin{equation}\label{hyp-Kisynski}
\forall t\in [0,T]~,~\forall u\in X^{1/2}~,~~ \langle 
H(t)u|u\rangle_{X} \geq \gamma \|u\|^2_{X^{1/2}} - \kappa \|u\|^2_{X}~.
\end{equation}
For any $u_0\in X^{1/2}$, there is a unique solution  $u \in \Cc^0([0,T], X^{1/2})\cap \Cc^1([0,T],X^{-1/2})$ of the equation 
\begin{equation}\label{eq-Kisynski}
i\partial_t u(t)=H(t)u(t)~~~~u(0)=u_0~.
\end{equation}
Moreover, $\|u_0\|_{X}=\|u(t)\|_{X}$ for all $t\in [0,t]$ and we may extend by 
density the flow of \eqref{eq-Kisynski} on $X$ as a unitary flow $U(t,s)$ such 
that $U(t,s)u(s)=u(t)$ for all solutions $u$ of \eqref{eq-Kisynski}.
If in addition $u_0\in D(H(0))$, then $u(t)$ belongs to $D(H(t))$ for all $t\in 
[0,T]$ and $u$ is of class $\Cc^1([0,T],X)$.
\end{theorem}

Let $\eta$, $I$ and $a$ be smooth controls. In our framework, the Hamiltonian 
\begin{equation}\label{def_H}
H(t):=-\ddd_{xx}^2 + I(t) \rho^{\eta(t)}(x-a(t)),\ \ \ \ \ \ t\in [0,T]
\end{equation}
is the time-dependent family of operators associated to the Schr\"odinger equation \eqref{eq} with $V$ defined in \eqref{eq_potentiel}. We notice that 
$$u\in H^1_0((0,1),\CC)\longmapsto V(t,\cdot)u \in H^{-1}((0,1),\CC)$$
is of class $\Cc^2$ with respect to $t$ due to the smoothness of $a$, $\eta$, $I$ and $\rho$. In addition to the positivity of each $H(t)$, this yields the hypothesis of Theorem \ref{th-Kisynski} (note that in a more singular setting of the free Schr\"odinger equation \eqref{eq_dmitry} with $\alpha$ going to $1$ and a moving cutting point $a$ the Cauchy problem would be more involved).  
\begin{coro}\label{coro_cauchy}
Let $T>0$ and let $V$ be defined in \eqref{eq_potentiel} with $a\in\Cc^2([0,T],(0,1))$, $I\in\Cc^2([0,T],\RR^+)$ and $\eta\in\Cc^2([0,T],\RR^+)$. For any $u_0\in H^1_0((0,1),\CC)$, Equation \eqref{eq} admits a unique solution $u \in \Cc^0([0,T],H^1_0((0,1),\CC))\cap \Cc^1([0,T],H^{-1}((0,1),\CC))$ such that 
$$\|u_0\|_{L^2}=\|u(t)\|_{L^2},\ \ \ \ \ \ \ \ \ \forall t\in [0,T].$$ 
The flow defined by \eqref{eq} can be unitary extended by density on $L^2((0,1),\CC)$. In addition, if $u_0\in H^2((0,1),\CC)\cap H^1_0((0,1),\CC)$, then $$u\in\Cc^0\Big([0,T],H^2((0,1),\CC)\cap H^1_0((0,1),\CC)\Big)  \cap \Cc^1([0,T],L^2((0,1),\CC)).$$
\end{coro}

We denote by $\G_s^t$ the unitary propagator in $L^2((0,1),\CC)$ generated by \eqref{eq}, as given by Corollary \ref{coro_cauchy}. The operator $\G_s^t$ represents the flow of \eqref{eq} in $[s,t]$ and, for any mild solution $u$ in $L^2((0,1),\CC)$ of the problem \eqref{eq}, we have
$\G_{s}^{t} u(s)=u(t)$.

\subsection{Adiabatic theory}\label{section_adiab}
In the following theorem, we present an important adiabatic result from Chapter IV of \cite{Bornemann} by using the notation adopted in Theorem \ref{th-Kisynski}. For further adiabatic results, we refer to the works \cite{Teufel,Schmid}. 

\begin{theorem}\label{th-Bornemann}
{\bf (Bornemann, 1998)}\\
Let the hypotheses of Theorem \ref{th-Kisynski} be satisfied for 
the times $t\in[0,1]$. Let $t\in[0,1]\mapsto\lambda(t)$ be a continuous curve such that $\lambda(t)$ for every $t\in [0,1]$ is in 
the discrete spectrum of $H(t)$, $\lambda(t)$ is a simple isolated eigenvalue 
for almost every $t\in[0,1]$ and there exists an associated family of orthogonal 
projections $P\in\Cc^1([0,1],\Lc(X))$ such that $$H(t)P(t)=\lambda(t)P(t),\ \ \ \  \ \ \ \forall t\in [0,1].$$
For any initial data $u_0\in X^{1/2}$ with $\|u_0\|_{X}=1$ 
and for any sequence $\epsilon\rightarrow 0$, the solutions 
$u_\epsilon \in \Cc^0([0,1],X^{1/2})\cap \Cc^1([0,1],X^{-1/2})$ of 
\begin{equation}\label{eq-Bornemann}
i\epsilon \partial_t 
u_\epsilon(t)=H(t)u_\epsilon(t)~~~~u_\epsilon(0)=u_0
\end{equation}
satisfy 
$$\langle P(1)u_{\epsilon}(1) |u_{\epsilon}(1)\rangle_{X} 
~~\xrightarrow[~~\epsilon\longrightarrow 0~~]{} \langle 
P(0)u_{\epsilon}(0) |u_{\epsilon}(0)\rangle_{X} $$
\end{theorem}

Theorem \ref{th-Bornemann} states the following fact. Let $t\mapsto\lambda(t)$ be a smooth curve of isolated simple eigenvalues and $t\mapsto\varphi(t)$ be a smooth curve of corresponding normalized eigenfunctions. If the initial data of \eqref{eq-Bornemann} is $u(0)=\varphi(0)$ and the dynamics induced by $H(t)$ with $t\in [0,1]$ is slow enough, then the final state $u(1)$ is as close as desired to the eigenfunction $\varphi(1)$. In our framework, we consider $H(t)$ defined by \eqref{def_H} where $\eta$, $I$ and $a$ are smooth controls. We notice that each $H(t)$ is positive with compact resolvent and there exists a Hilbert basis of $L^2((0,1),\CC)$ made by its eigenfunctions. In order to apply Theorem \ref{th-Bornemann} to the Schr\"odinger equation \eqref{eq} with $V$ defined in \eqref{eq_potentiel}, we ensure the simplicity of the eigenvalues of $H(t)$ at each fixed $t\in[0,1]$, as given by the following proposition. 
\begin{prop}\label{prop-simple}
For any $a\in (0,1)$, $I\geq 0$ and $\eta>0$, the eigenvalues of $-\ddd_{xx}^2 + I \rho^{\eta}(x-a)$ are simple. 
\end{prop}
\begin{demo}
If $\phi$ and $\psi\in H^2((0,1),\CC)\cap H^1_0((0,1),\CC)$ are two eigenfunctions corresponding to 
the same eigenvalue $\lambda$, then they both satisfy the same second order differential equation 
$$f''(x)=I\rho^\eta(x-a)f(x)+\lambda f(x)~.$$
Since $\rho$ is continuous, they are both of class $\Cc^2$ and thus classical solutions of the above ODE. Moreover, $\phi(0)=\psi(0)=0$ and there exist $\alpha,\beta\in\CC$ such that $\alpha\phi'(0)=\beta \psi'(0)$. Thus, $\alpha\phi$ and $\beta\psi$ satisfy the same second order ODE with the same initial data and $\alpha\phi\equiv\beta\psi$\qedhere. 
\end{demo}

As a consequence of the result of Proposition \ref{prop-simple}, we can respectively define by $$\big(\lambda_k(t)\big)_{k\in\NN^*},\ \ \ \  \ \ \ \ \big(\phi_k(t)\big)_{k\in\NN^*}$$
the ordered sequence of eigenvalues of $H(t)=-\ddd_{xx}^2 + I(t) \rho^{\eta(t)}(x-a(t))$ for every $t\in[0,1]$ and a Hilbert basis made by corresponding eigenfunctions. In the next proposition, we ensure the validity of the adiabatic result of Theorem \ref{th-Bornemann} for the Schr\"odinger equation \eqref{eq} with $V$ defined in \eqref{eq_potentiel}.

\begin{prop}\label{prop-adiab}
Let $N\in\NN^*$. For every $\varepsilon>0$, $\widetilde\eta\in\Cc^2([0,1],\RR^+)$, $\widetilde I\in\Cc^2([0,1],\RR^+)$ and $\widetilde a\in\Cc^2([0,1],(0,1))$, there exists $T^*>0$ such that, for every $T\geq T^*$, the following property holds. Take $\eta(t)=\widetilde\eta\big(\frac{t}{T}\big)$, $I(t)=\widetilde I\big(\frac{t}{T}\big)$ and $a(t)=\widetilde a\big(\frac{t}{T}\big)$, then the flow of the corresponding Schr\"odinger equation \eqref{eq} satisfies
$$\forall\, k\leq N~,~~\exists\, \alpha_k\in\CC~:~~|\alpha_k|=1~~\text{ and }~~
\big\|\Gamma_0^T\phi_{k}(0)-\alpha_k\phi_{k}(T)\big\|_{L^2}\leq \varepsilon~.$$
\end{prop}
\begin{demo}
The family of self-adjoint operators $H(t)$ is smooth in $t\in [0,1]$ and its eigenvalues are simple. By classical spectral theory, the eigenvalues $\big(\lambda_k(t)\big)_{k\in\NN^*}$ with $t\in[0,1]$ form smooth curves and the basis made by corresponding eigenfunctions $\big(\phi_k(t)\big)_{k\in\NN^*}$ can be chosen smooth with respect to $t\in[0,1]$, up to adjusting the phases (see \cite{Kato}). It remains to apply Theorem \ref{th-Bornemann} with respect to the chosen Hamiltonian and to notice that if $u_\epsilon$ solves \eqref{eq-Bornemann}, then $u(t)=u_\epsilon(\epsilon t)$ solves our main equation \eqref{eq} with $a(t)=\tilde a(\epsilon t)$, $I(t)=\tilde I(\epsilon t)$ and $\eta(t)=\tilde\eta(\epsilon t)$. To conclude, we choose $\epsilon$ so that Theorem \ref{th-Bornemann} yields an error at most $\varepsilon$ for the first $N$ eigenmodes and we set $T^*=1/\epsilon$.
\end{demo}


\section{Spectral analysis for large $I$ and $\eta$}\label{section_spectre}

To mimic the dynamics of the model \eqref{eq_dmitry}, we need to have a very high and very sharp potential wall. In other words, in proving Theorem \ref{th}, we consider both $I$ and $\eta$ very large. 

\vspace{1mm}

In this section, we study the spectrum of the positive self-adjoint Hamiltonian
$$H^{I,\eta,a}=-\ddd_{xx}^2 + I \rho^{\eta}(x-a),\ \ \ \ \text{ with }\ \ \ \ \eta>0,\ I>0,\ a\in (0,1).$$ 
Heuristically speaking, we expect that, from a spectral point of view, the limit of $H^{I,\eta,a}$ for 
$I,\eta\rightarrow +\infty$ is the operator $H^{\infty,a}$ defined as follows
$$H^{\infty, a}=-\ddd_{xx}^2,\ \ \ \ \ \ \ \ \ \ \ D(H^{\infty, a})=H^2\Big((0,a)\cup (a,1),\CC\Big) \cap  H^1_0\Big((0,a)\cup (a,1),\CC\Big).$$
Such operator corresponds to the Laplacian on a split segment $(0,a)\cup(a,1)$ with Dirichlet boundary conditions as considered in \cite{Dmitry}.
We respectively denote by 
$$\big(\lambda_k^{I,\eta,a}\big)_{k\in\NN^*},\ \ \ \ \  \ \ \ \ \big(\phi_k^{I,\eta,a}\big)_{k\in\NN^*}$$ 
the ordered sequence of eigenvalues of $H^{I,\eta,a}$ and a Hilbert basis of $L^2((0,1),\CC)$ made by corresponding eigenfunctions. 
The spectrum of $H^{\infty, a}$ is composed by the numbers $(\mu_p^l(a))_{p\in\NN^*}$ and $(\mu_q^r(a))_{q\in\NN^*}$ defined in \eqref{def_mu} 
and we consider a Hilbert basis made by corresponding eigenfunctions given by
\begin{align}
\varphi_p^l(x)&=\sqrt{\frac{2}{a}}\sin\Big(\frac{p\pi}{a} x \Big)\Un_{x\in[0,a]},\label{def_varphi_l}\\
\varphi_q^r(x)&= \sqrt{\frac{2}{1-a}}\sin\Big(\frac{q\pi}{1-a}(1-x)\Big)\Un_{x\in[a,1]}. \label{def_varphi_r}
\end{align}
where $\Un_{x\in J}$ is equal to $1$ in $J\subset (0,1)$ and to $0$ in $(0,1)\setminus J$. We denote as $$(\lambda^{\infty,a}_k)_{k\in\NN^*},\ \ \ \mbox{ and  } \ \ \ \  \ \ \ \  (\phi^{\infty,a}_k)_{k\in\NN^*}$$ the ordered spectrum of $H^{\infty, a}$ obtained by reordering the eigenvalues given by \eqref{def_mu} in ascending order and,
respectively, a Hilbert basis of $L^2((0,1),\CC)$ made by corresponding eigenfunctions (defined as \eqref{def_varphi_l} or \eqref{def_varphi_r}). We notice that $H^{\infty, a}$ may have multiple eigenvalues, unlike the operator $H^{I,\eta,a}$ (see Proposition \ref{prop-simple}). More precisely, $H^{\infty, a}$ has at most double eigenvalues appearing if and only if $a$ is rational.  

\vspace{3mm}

The purpose of this section is to prove the following result which shows that the spectrum 
$(\lambda_k^{I,\eta,a})_{k\in\NN^*}$ of $H^{I,\eta,a}$ converges to the spectrum of $H^{\infty, a}$ when $\eta$ and $I$ go to $+\infty$. We will only consider the case where $I$ and $\eta$ are of the same order, in other words when there exists $\delta>0$ such that  
\begin{equation}\label{Ieta}
\delta \eta ~\leq~I~\leq~\frac 1\delta \eta~. 
\end{equation}
\begin{theorem}\label{th-cv-spectre} 
For each $k\in\NN^*$ and $\delta>0$, there exists a constant $C_{k,\delta}>0$ (which depends continuously on $a\in(0,1)$ and on the shape of $\rho$) such that, for all
$\eta \geq 1$ and $I\geq 1$ satisfying \eqref{Ieta}, we have
$$\lambda^{\infty,a}_k\,-\,\frac{C_{k,\delta}}{\sqrt{\eta}}~\leq~
\lambda_k^{I,\eta,a}~\leq~\lambda^{\infty,a}_k + \frac{C_{k,\delta}}{\eta}~.$$
If $\lambda_k^{\infty,a}$ is a simple eigenvalue of $H^{\infty, a}$ with a normalized eigenfunction $\phi^{\infty,a}_k$, then there exists $\alpha_k^{I,\eta,a}\in\CC$ such that $|\alpha_k^{\eta,a}|=1$ and
$$\|\phi_k^{I,\eta,a} - \alpha_k^{I,\eta,a} \phi^{\infty,a}_k\|_{L^2}~\leq~\frac {C_{k,\delta}}{\sqrt{\eta}}~. 
$$
If $\lambda_k^{\infty,a}=\lambda_{k+1}^{\infty,a}$ is a double eigenvalue of $H^{\infty, a}$ corresponding to a pair of normalized eigenfunction $\phi^{\infty,a}_k$ and $\phi^{\infty,a}_{k+1}$, then there exists $\alpha_k^{I,\eta,a}\in\CC$ and $\beta_k^{I,\eta,a}\in\CC$ such that 
$|\alpha_k^{I\eta,a}|^2 + |\beta_k^{I,\eta,a}|^2 =1$ and
$$\|\phi_k^{I,\eta,a} - \alpha_k^{I,\eta,a} \phi^{\infty,a}_k - \beta_k^{I,\eta,a} \phi^{\infty,a}_{k+1}\|_{L^2}~\leq~\frac {C_{k,\delta}}{\sqrt{\eta}}~. 
$$
\end{theorem}

\noindent We split the proof of Theorem \ref{th-cv-spectre} into several lemmas given below.
\begin{lemma}\label{lemma-spectre-1}
For each $k\in\NN^*$, $a\in (0,1)$, $I>0$ and $\eta >0$, we have that
$$\lambda_k^{I,\eta,a}~\leq~\lambda^{\infty,a}_k \left(1+ 2k \frac{I}{\min(a,1-a)\eta^2}\right)~.$$
\end{lemma}
\begin{demo}
The main tool of the proof is the min-max 
theorem (see \cite[Theorem\ XIII.1]{Reed-Simon}). We introduce the Rayleigh quotient 
$$R^{I,\eta,a}(u)~=~ \frac{\langle 
H^{I,\eta,a} u |u\rangle_{L^2}}{\|u\|^2_{L^2}}~=~ \frac{\int_0^1\big( 
|\ddd_xu(x)|^2 + I \rho^\eta(x-a)|u(x)|^2\big)\dd x}{\int_0^1 |u(x)|^2\dd x} ~.$$
The min-max principle yields that, for every $k\in\NN^*$,
$$\lambda_k^{I,\eta,a} ~=~\min \big\{ \max_{u\in E} R^{I,\eta,a}(u)~|~E\subset H^1_0,\ dim(E)=k \big\}~.$$ 

Thus, 
$$\lambda_k^{I,\eta,a} ~\leq \max_{u\in E} R^{I,\eta,a}(u),$$ 
where $E=\text{span}(\{\phi^{\infty,a}_j, j=1\ldots k\})$ is the space spanned by the first $k$ eigenfunctions of the limit operator $H^{\infty,a}$. First, we estimate $R^{I,\eta,a}(\phi^{\infty,a}_j)$.
Let $\phi^{\infty,a}_j=\sin \big(\frac{p\pi x}{a}\big)\Un_{x\in(0,a)}$ with $p\in\NN^*$ (a similar computation is valid for the eigenfunctions of the type 
$\sin \big(\frac{q\pi}{1-a} (1-x)\big)\Un_{x\in(a,1)}$ with $q\in\NN^*$). 
Since $\big|\sin\big(\frac{p\pi x}{a}\big)\big|\leq \frac{p\pi}{a\eta}$ in $[-\frac{1}{\eta},0]$ and $\int_{-\frac{1}{\eta}}^{0}\rho^{\eta}(x)\dd x=\int_{-1}^{0}\rho(x)\dd x\leq 1$, we have 
$$ \int_0^1 I \rho^\eta(x-a)|\phi^{\infty,a}_j(x)|^2 \dd x = I 
\int_{-\frac{1}{\eta}}^{0} \sin^2\big(\frac{p\pi x}{a}\big) \rho^{\eta}(x)\dd x
\leq I \frac{p^2\pi^2}{a^2\eta^2} = I \frac{\lambda^{\infty,a}_j}{\eta^2} \leq 
I \frac{\lambda^{\infty,a}_k}{\eta^2}.
$$
We also have
$$ \int_0^1 |\partial_x \phi^{\infty,a}_j(x)|^2 \dd x=
\lambda^{\infty,a}_j \int_0^1 |\phi^{\infty,a}_j(x)|^2\dd x= \frac {\gamma_j}2 \lambda^{\infty,a}_j$$ 
where $\gamma_j$ is either $a$ or $1-a$, depending of the type of the eigenfunction. Moreover, notice that the functions $\phi^{\infty,a}_j$ are orthogonal, both for $L^2$ and $H^1$ scalar products. Thus, we obtain that, for any $u=\sum_{j=1}^k c_j\phi^{\infty,a}_j \in E$,
\begin{align*}
R^{I,\eta,a}(u)&=\frac{\frac 12 \sum |c_j|^2{\gamma_j}\lambda^{\infty,a}_j + 
\int_0^1 I \rho^{\eta}(x-a)|\sum c_j\phi^{\infty,a}_j(x)|^2 \dd x} { \frac 12 \sum_j |c_j|^2\gamma_j}\\
&\leq \lambda^{\infty,a}_k + \frac{2}{\sum_j |c_j|^2\gamma_j} \int_0^1 I \rho^{\eta}(x-a) \;k 
\sum |c_j\phi^{\infty,a}_j(x)|^2 \dd x\\
&\leq \lambda^{\infty,a}_k + 
\frac{2k}{\sum_j |c_j|^2\gamma_j} \sum |c_j|^2 \int_0^1 I \rho^{\eta}(x-a)|\phi^{\infty,a}_j(x)|^2 \dd x\\
&\leq \lambda^{\infty,a}_k + 
\frac{2k}{\sum_j |c_j|^2\gamma_j} \sum |c_j|^2 I \frac{\lambda^{\infty,a}_k}{\eta^2}\\
& \leq \lambda^{\infty,a}_k \left(1+ 2k \frac{I}{\min(a,1-a)\eta^2}\right)
\end{align*}
Then, the min-max principle yields the result.\ \ \qedhere
\end{demo}

\begin{lemma}\label{lemma-spectre-2}
For any $k\in\NN^*$, we have that 
$$\lambda_k^{I,\eta,a}~\xrightarrow[~~I,\eta\longrightarrow+\infty~~]{}\lambda^{\infty,a}_k~,$$
provided $I=o(\eta^2)$. Moreover, for any sequences $(I_n)_{n\in\NN^*}$ and $(\eta_n)_{n\in\NN^*}$ going to $+\infty$ such that $I_n=o(\eta_n^2)$, there exist two 
subsequences $(I_{n_j})_{j\in\NN^*}$, $(\eta_{n_j})_{j\in\NN^*}$ and a normalized eigenfunction 
$\phi^{\infty,a}_k$ of $H^{\infty, a}$ for the eigenvalue $\lambda^{\infty,a}_k$ such that 
$$\phi_k^{I_{n_j},\eta_{n_j},a}~\xrightarrow[j\longrightarrow+\infty]{}
\phi^{\infty,a}_k$$ 
 weakly in $H^1_0((0,1),\CC)$ and strongly in $H^{1-\varepsilon}((0,1),\CC)$ with $\varepsilon\in(0,1]$.
\end{lemma}
\begin{demo}
Let $(I_n,\eta_n)_{n\in\NN^*}$ be a sequence of parameters going to $+\infty$, such that 
$I_n = o(\eta_n^2)$. Extracting a subsequence, if necessary, we can assume that 
$\lambda_k^{I_n,\eta_n,a}\xrightarrow[n\rightarrow+\infty]{}r_k\in
[0,\lambda^{\infty,a}_k]$, thanks to the bound of Lemma \ref{lemma-spectre-1}. 
Now, we use the identity $\la\phi_k^{I,\eta,a}|H^{\eta,a} 
\phi_k^{I,\eta,a}\ra_{L^2}=\lambda_k^{I,\eta,a}$, which leads to
\begin{equation}\label{eq-proof-th-cv-spectre}
\int_0^1 |\ddd_x 
\phi_k^{I,\eta,a}(x)|^2 \dd x~+
~I \int_0^1 \rho^{\eta}(x)|\phi_k^{I,\eta,a}(x)|^2 \dd x~
=~\lambda_k^{I,\eta,a}~\leq~C_k
\end{equation}
for some constant $C_k>0$ provided by Lemma \eqref{lemma-spectre-1} and $I=\Oc(\eta^2)$.

Thus, $\|\phi^{I,\eta,a}_k\|^2_{H^1_0} \leq C_k$ and, by compactness of the weak topology, we can assume that the sequence $\phi_k^{I_n,\eta_n,a}$ weakly converges in $H^1((0,1),\CC)$ to a function $\psi_k$ (up to extracting a subsequence). By compactness of Sobolev embeddings, we can also assume that the convergence is strong in $H^{1-\varepsilon}((0,1),\CC)$ for every $\varepsilon>0$. 
In particular, the limit functions $(\psi_k)_{k\in\NN^*}$ form an orthonormal family of $L^2((0,1),\CC)$. 

Now, we claim that the functions $\phi_k^{I,\eta,a}$ must be small in a neighborhood of $a$. Indeed, 
$\|\phi_k^{I,\eta,a}\|^2_{H^1_0} \leq C_k$ yields that the functions $(\phi_k^{I,\eta,a})_{k\in\NN^*}$ are uniformly $\frac12-$H\"older continuous. Thus, 
\begin{align*}
|\phi_k^{I,\eta,a}(a)|^2 & = |\phi_k^{I,\eta,a}(x) + (\phi_k^{I,\eta,a}(a)-\phi_k^{I,\eta,a}(x))|^2\\
&\leq 2 |\phi_k^{I,\eta,a}(x)|^2 + 2 |\phi_k^{I,\eta,a}(x)-\phi_k^{I,\eta,a}(a)|^2\leq 
2 |\phi_k^{I,\eta,a}(x)|^2 + \Oc(\eta^{-1})
\end{align*}
for any $x\in [a-1/\eta,a+1/\eta]$ (i.e., in the support of $\rho^\eta$). This implies that
$$|\phi_k^{I,\eta,a}(a)|^2 \leq \int_0^1 \rho^{\eta}(x)|\phi_k^{I,\eta,a}(x)|^2 \dd x + \Oc(\eta^{-1})$$
Thus, the bound \eqref{eq-proof-th-cv-spectre} gives
$$|\phi_k^{I_n,\eta_n,a}(a)|\leq \frac{C}{\min\{\sqrt{\eta_n},\sqrt{I_n}\}}$$
and, by the $\frac12-$H\"older continuity, the same bound holds for all $x\in [a-1/\eta,a+1/\eta]$;
the constant $C$ depends on $k$, $a$, and the choice of $\rho$.
The strong convergence in $H^{\frac{3}{4}}((0,1),\CC)$ and in $\Cc^0((0,1),\CC)$ implies that  $\psi_k(a)=0$ for all $k\in\NN^*$. 

Now, let $\varphi\in\Cc^\infty_0((0,a)\cup(a,1),\CC)$ be a test function. For $n\in\NN^*$ large, we have $\varphi\rho^{\eta_n}\equiv 0$ and
\begin{align*}
\int_0^1 \grad \phi_k^{I_n,\eta_n,a}(x)\grad \varphi(x)\dd x &=
\la (-\partial_{xx}^2)\phi_k^{I_n,\eta_n,a}|\varphi\ra_{L^2}=
\la H^{I_n,\eta_n,a}\phi_k^{I_n,\eta_n,a}|\varphi\ra_{L^2}\\
&= \lambda_k^{I_n,\eta_n,a} \int_0^1 \phi_k^{I_n,\eta_n,a}(x)\varphi(x)\dd x.
\end{align*}
Passing to the weak $H^1-$limit, for every $\varphi\in\Cc^\infty_0((0,a)\cup(a,1),\CC)$, we find
$$\int_0^1 \grad \psi_k(x) \grad 
\varphi(x) \dd x = r_k \int_0^1 \psi_k(x)\varphi(x)\dd x~.$$
The above equality and the fact that $\psi_k(a)=0$ show that $r_k$ (the limit of a subsequence
$\big(\lambda^{I_{n_j},\eta_{n_j},a}_k\big)_{j\in\NN^*}$) is, for any choice of the subsequence $(n_j)_{j\in\NN^*}$, an eigenvalue of $H^{\infty, a}$ and $\psi_k$ is a corresponding eigenfunction. Moreover, we have 
$r_s\leq \lambda^{\infty,a}_k$ for all $s\leq k$ and, by the linear independence of the 
limit functions $(\psi_k)_{k\in\NN^*}$, we obtain that the numbers $(r_s)_{1\leq s \leq k}$ must 
be the first $k$ eigenvalues of $H^{\infty, a}$, implying that $r_k$ is the $k$-th eigenvalue and $\psi_k$
is the $k$-th eigenfunction $\phi^{\infty,a}_k$.
\end{demo}

\begin{lemma}\label{lemma-spectre-3}
For each $k\in\NN^*$ and $\delta>0$, there exists a constant $C_{k,\delta}>0$ (depending continuously on $a\in(0,1)$, on the shape of $\rho$) such that, for all $\eta\geq 1$ and $I\geq 1$ satisfying \eqref{Ieta},
$$\lambda_k^{I,\eta,a} ~\geq~\lambda_k^{\infty,a}\,-\,\frac{C_{k,\delta}}{\sqrt{\eta}}.$$
\end{lemma}
\begin{demo}
From the proof of Lemma \ref{lemma-spectre-2}, we know that there exists
${C_{k,\delta}}>0$ such that  
\begin{equation}\label{eq-proof-th-cv-spectre-2}
\big|\phi_k^{I,\eta,a}(x)\big|~\leq~\frac {C_{k,\delta}}{\sqrt{\eta}}~~~\text{ for all }
x\in \big[a-1/\eta,a+1/\eta \big]~.
\end{equation}
In particular, the $L^2-$norm of $\phi_k^{I,\eta,a}$ in $\big[a-1/\eta,a+1/\eta \big]$ goes to zero as 
$\eta$ and $I$ grow. As a consequence, the $L^2-$norm of $\phi_k^{I,\eta,a}$ is at least almost half supported in one of the 
intervals $\big[0,a-{1}/{\eta}\big]$ or $\big[a+1/{\eta},1\big]$. Without loss of generality, we assume 
that $\|\phi_k^{I,\eta,a}\|_{L^2((0,a-\frac{1}{\eta}),\CC)}\geq \frac{1}{4}$. Since $\ddd_{xx}^2 
\phi_k^{I,\eta,a} = -\lambda_k^{I,\eta,a} \phi_k^{I,\eta,a}$ in $\big[0,a-{1}/{\eta}\big]$ and due to the 
Dirichlet boundary condition at $x=0$, we have that  
$\phi_k^{I,\eta,a}(x)=A_k \sin(\sqrt{\lambda_k^{I,\eta,a}}x)$ when $x\in \big[0,a-\frac{1}{\eta}\big]$
with a constant $A_k\in\CC$. Since the $L^2$-norm of $\phi_k^{I,\eta,a}$ on the interval 
$\big[0,a-\frac{1}{\eta}\big]$ is bounded away from zero, it follows that $A_k$ stays bounded away from zero as $\eta$ and $I$ grow.
Thus, from 
\eqref{eq-proof-th-cv-spectre-2}, we must have $\sqrt{\lambda_k^{I,\eta,a}}=\frac{k'\pi}{a} 
+ \Oc\big(1/{\sqrt{\eta}}\big)$, and $k'=k$ for $\eta$ large enough thanks to the convergence obtained in Lemma 
\ref{lemma-spectre-2}. \end{demo}

\begin{lemma}\label{lemma-spectre-4}
Let $k\in\NN^*$ and $\delta>0$. For each $\eta\geq 1$ and $I\geq 1$ satisfying \eqref{Ieta}, there exists a normalized eigenfunction 
$\widetilde\phi^{\infty,a}_k$ of $H^{\infty, a}$ corresponding to the eigenvalue $\lambda^{\infty,a}_k$ such that 
\begin{align}\label{lemma-spectre-4-eq-1}\|\phi_k^{I,\eta,a} - \widetilde\phi^{\infty,a}_k\|_{L^2}~\leq~\frac {C_{k,\delta}}{\sqrt{\eta}}~\end{align}
where $C_{k,\delta}>0$ depends on $k$ and $\delta$ and also on $a\in(0,1)$ and the shape of $\rho$.
\end{lemma}
\begin{demo}
Let $k\in\NN^*$ and assume that \eqref{Ieta} holds. From the proof of the previous lemmas, we know that 
\begin{equation}\label{lemma-spectre-4-eq-2}
 \sqrt{\lambda_k^{I,\eta,a}}~=~\sqrt{\lambda^{\infty,a}_k} ~+~\Oc\left( \frac{1}{\sqrt{\eta}}\right),
\end{equation}
Since $\ddd_{xx}^2  \phi_k^{I,\eta,a} = -\lambda_k^{I,\eta,a} \phi_k^{I,\eta,a}$ outside of $\big[a-{1}/{\eta}, a+{1}/{\eta}\big]$ and due to the Dirichlet boundary conditions at $x=0$ and $x=1$,
we have that
there exist $A_k,B_k\in\CC$ such that 
\begin{align}
\phi_k^{I,\eta,a}(x)~=~&A_k 
\sin\Big(\sqrt{\lambda_k^{I,\eta,a}}x\Big)\Un_{x\in [0,a-1/{\eta}]}~+~
\Oc\big(\frac 1{\sqrt \eta}\big)\Un_{x\in [a-{1}/{\eta},a+ 1/\eta]} \nonumber \\
&~+~B_k \sin\Big(\sqrt{\lambda_k^{I,\eta,a}}(1-x)\Big)\Un_{x\in[a+{1}/{\eta},1]}\label{lemma-spectre-4-eq-3} 
\end{align}
(the estimate on $\phi_k^{I,\eta,a}$ at $x \in \big[a-{1}/{\eta}, a+{1}/{\eta}\big]$ is given by 
(\ref{eq-proof-th-cv-spectre-2})).

By normalization, $a |A_k|^2+ (1-a)|B_k|^2=2+\Oc\big(\frac{1}{\sqrt\eta}\big)$. Now, we have two possibilities. If $\lambda^{\infty,a}_k=\lambda^{\infty,a}_{k+1}$ is a double eigenvalue, then we take
$$\widetilde\phi^{\infty,a}_k=\alpha \sin\big(\sqrt{\lambda^{\infty,a}_k } x \big) \Un_{x\in[0,a]} + \beta\sin\big(\sqrt{\lambda^{\infty,a}_k}(1-x) \big)\Un_{x\in[a,1]}$$
with 
$$ \alpha=\frac{\sqrt{2}A_k}{\sqrt{a|A_k|^2+(1-a)|B_k|^2}}~~~\text{ and }~~~\beta=\frac{\sqrt{2}B_k}{\sqrt{a|A_k|^2+(1-a)|B_k|^2}}.$$
Obviously, this is the sought normalized eigenfunction of $H^{\infty, a}$ satisfying 
\eqref{lemma-spectre-4-eq-1} due to \eqref{lemma-spectre-4-eq-2} and \eqref{lemma-spectre-4-eq-3}. 

In the case where $\lambda^{\infty,a}_k$ is a simple eigenvalue, one of the functions 
$\sin(\sqrt{\lambda^{\infty,a}_k} x)$ or $\sin(\sqrt{\lambda^{\infty,a}_k}(1-x))$ does not satisfy 
the Dirichlet condition at $x=a$. Without loss of generality, we assume that 
$$\sin\Big(\sqrt{\lambda^{\infty,a}_k} a\Big)=0\ \ \ \text{ and }\ \ \ \ \sin\Big(\sqrt{\lambda^{\infty,a}_k}(1-a)\Big)\neq 0.$$  
By \eqref{eq-proof-th-cv-spectre-2}, we conclude that 
$B_k^{\eta,a}=\Oc\big(\frac{1}{\sqrt{\eta}}\big)$. Therefore, \eqref{lemma-spectre-4-eq-1} holds for
$$\widetilde\phi^{\infty,a}_k=\frac{A_k}{|A_k|} \sqrt{\frac 2a}\sin\big(\sqrt{\lambda^{\infty,a}_k} x \big) 
\Un_{x\in[0,a]}~=\frac{A_k}{|A_k|} \phi^{\infty,a}_k.$$
\end{demo}


\section{Short-time movement of the high potential wall}\label{section_fast_translation}

In this section, we consider $V(t,x)=I(t) \rho^{\eta(t)}(x-a(t))$ with $I$, $\eta$ and $a$ being smooth controls. For $\eta$ large, we can see $V$ as a very thin (of width $\eta^{-1}$) moving wall. When the parameter $I$ that determines the height of the wall is large the norm of the map
$$u\in H^1_0((0,1),\CC)~\longmapsto~ V(x,t)u=I(t)\rho^{\eta(t)}(x-a(t))u \in H^{-1}((0,1),\CC)$$ 
is large. Therefore, even a small displacement of $a(t)$ may, in principle, cause
a strong modification of the initial state. However, the following lemma shows that if the initial condition is  localized outside the support interval $\big[a(t)-\frac{1}{\eta},a(t)+\frac{1}{\eta}\big]$ of the potential $V$,
then the speed with which the solution deviates from the initial condition is uniformly bounded for all $\eta$ and $I$.
\begin{lemma}\label{lemma-compare}
Let $\eta(t)\geq\eta_*\in\RR^+$ for all $t$ from some time interval $[t_1,t_2]$, and let $u(t)$ be the solution of the Schr\"odinger equation \eqref{eq} in this interval with $L^2((0,1),\CC)$ initial data. Take any
function 
$$\psi\in C^1([t_1,t_2],H^2((0,1),\CC)\cap H^1_0((0,1),\CC))$$ 
which, for every $t\in[t_1,t_2]$, vanishes in $\big[a(t)-\frac{1}{\eta_*},a(t)+\frac{1}{\eta_*}\big]$
(so $\psi V \equiv 0$). Then, for all $t\in [t_1,t_2]$,
\begin{equation}\label{eq-bootstrap-intermediate}
\|u(t)-\psi(t)\|_{L^2}^2 \leq \|u(t_1)-\psi(t_1)\|_{L^2}^2 + C(t_2-t_1),
\end{equation}
with $C$ independent of the choice of functions $\eta$, $I$, $a$ and given by
$$C=2\sup_{t\in [t_1,t_2]} \big(\,\|u(t_1)\|_{L^2} \|\partial_{xx}^2\psi(t)\|_{L^2} + 
(\|\psi(t)\|_{L^2}+\|u(t_1)\|_{L^2})\|\partial_t \psi(t)\|_{L^2}\,\big)~.$$
\end{lemma}
\begin{demo}
As the evolution operator of the Schr\"odinger equation is unitary, $\|u(t)\|_{L^2}$ is constant in time. Let us, first, assume that $u(t_1)\in D(-\ddd_{xx}^2+V(t_1,\cdot))$, which implies regularity of the solution with respect to time as shown in Corollary \ref{coro_cauchy}. 
As $V(t,\cdot)\psi(t,\cdot)\equiv 0$ for every $t\in [t_1,t_2]$, we have 
\begin{align*}
&\partial_t \frac 12 \|u(t)-\psi(t) \|^2_{L^2}= \partial_t \left(\frac 12\|\psi(t)\|^2  
-\Re\big( \langle u(t) |  \psi(t) \rangle_{L^2}  \big)\right)\\
& \leq \|\psi(t)\|_{L^2}\|\partial_t\psi(t)\|_{L^2} + |\langle (-\ddd_{xx}^2+V(t,\cdot))u(t) | \psi(t)\rangle_{L^2}| 
+ \|u(t)\|_{L^2}\|\partial_t \psi(t)\|_{L^2}\\
&\leq |\langle \partial_{xx}^2u(t) | \psi(t)\rangle_{L^2}|+ (\|\psi(t)\|_{L^2}+\|u(t_1)\|_{L^2})\|\partial_t \psi(t)\|_{L^2}\\
& \leq |\langle u(t) | \partial_{xx}^2\psi(t)\rangle_{L^2}| + 
(\|\psi(t)\|_{L^2}+\|u(t_1)\|_{L^2})\|\partial_t \psi(t)\|_{L^2}\\
&\leq \|u(t_1)\|_{L^2} \|\partial_{xx}^2\psi(t)\|_{L^2} + (\|\psi(t)\|_{L^2}+\|u(t_1)\|_{L^2})\|\partial_t \psi(t)\|_{L^2}~.
\end{align*}
Now the validity of estimate \eqref{eq-bootstrap-intermediate} follows from that, 
for every $t\in [t_1,t_2]$,
\begin{equation*}\|u(t)-\psi(t)\|_{L^2}^2 \leq \|u(t_1)-\psi(t_1)\|_{L^2}^2 + |t_2-t_1| \sup_{t\in [t_1,t_2]}\big|\partial_t \|u(t)-\psi(t) \|^2_{L^2}\big|. \end{equation*}
We conclude the proof by noticing that the estimate is extended by density 
to any $u(t_1)\in L^2((0,1),\CC)$.
\end{demo}

Recall that we denote the eigenvalues of the Laplacian operator on the split interval $(0,a(t))\cup(a(t),1)$
for any frozen value of $t$ as $\lambda^{\infty,a(t)}_k$; for each $k\in\NN^*$, we have $\lambda^{\infty,a}_k\in\{\mu_p^l(a)\}_{p\in\NN^*}$ or $\lambda^{\infty,a}_k\in\{\mu_q^r(a)\}_{q\in\NN^*}$ 
where $\mu_p^l$ and $\mu_q^r$ are introduced in \eqref{def_mu}. The corresponding
eigenfunctions are denoted as $\phi_k^{\infty,a(t)}$ and are equal, respectively, to $\varphi_p^l$ or $\varphi_q^r$ defined by \eqref{def_varphi_l} or \eqref{def_varphi_r}. The functions 
$\varphi_p^l$ and $\varphi_q^r$ evolve smoothly with $a$ and are mostly localized outside the support
interval of $V$, so if we consider any of these functions as an initial condition for the equation \eqref{eq},
then we can extract from Lemma \ref{lemma-compare} (see Proposition below) that the solution will not deviate far from it on a sufficiently short interval of time, uniformly for all large $\eta$ and all $I$. 

\begin{figure}[H]
\begin{center}
\resizebox{0.8\textwidth}{80mm}{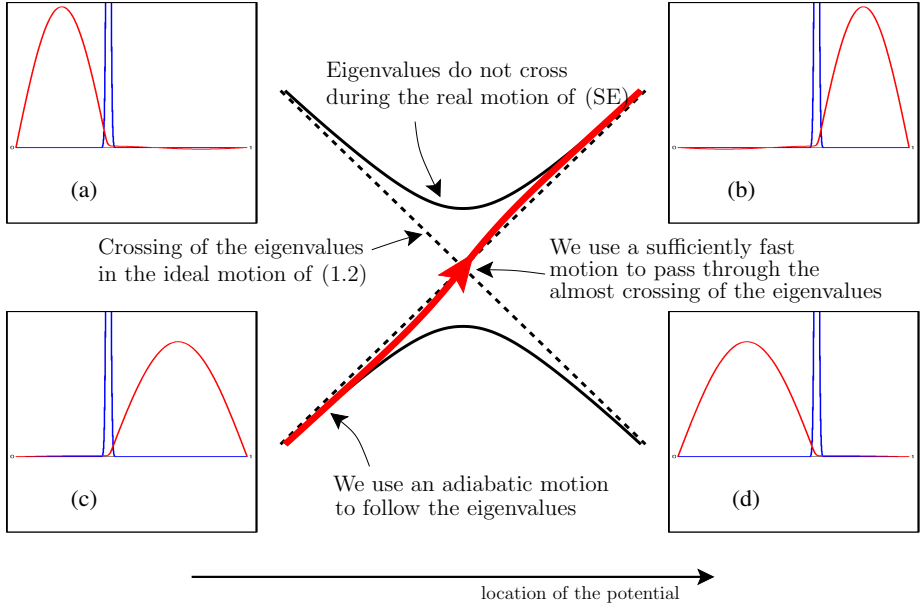}
\end{center}
\caption{\it The evolution of the eigenmode {(c)} during the motion of the potential wall. On one hand, when the wall moves very slow and the time it takes to pass the crossing of the eigenvalues of the split Laplacian is large, the motion is adiabatic and the state is deformed as shown in {(d)} - it tunnels from one  side of the wall $V$ to the other side (because the spectral curves of $-\ddd_{xx}^2+V$ do not intersect). On the other hand, when we translate the wall in a short time, it is possible to ``jump'' from the lower spectral curve to the higher as showed by the red line in the figure. In such case, the state is deformed 
as in {(b)}, because this state is close to the mode shown in {(c)}.
}\label{fig_crossing} 
\end{figure}

In particular, if $a_*\in (0,1)$ is a crossing point such that $\lambda_k^{\infty,a_*}=\lambda_{k+1}^{\infty,a_*}$ for some $k\in\NN^*$, then for any
small $\delta$, the split Laplacian eigenfunctions  $\phi_k^{\infty,a_*-\delta}$ and $\phi_{k+1}^{\infty,a_*+\delta}$ correspond to the same mode $\varphi_p^l$ or $\varphi_q^r$.
Therefore, if $a(t_1)=a_*-\delta$ and $a(t_2)=a_*+\delta$, then the flow of \eqref{eq} on
the interval $[t_1,t_2]$ will take the function $\phi_k^{\infty,a(t_1)}$ close to $\phi_{k+1}^{\infty,a(t_2)}$
(and the function $\phi_{k+1}^{\infty,a(t_1)}$ close to $\phi_{k}^{\infty,a(t_2)}$), provided $\delta$ and $(t_2-t_1)$ are small enough (see Fig. \ref{fig_crossing}).

This gives exactly the permutation of eigenmodes introduced in Theorem \ref{th}. The precise statement is given by the following 
\begin{prop}\label{prop-compare-eigenfunctions1} Let $a_*\in (0,1)\cap\QQ$.
Take any $\varepsilon>0$, $\kappa>0$, and $N\in\NN^*$ such that $\lambda_{N+1}^{a_*}\neq\lambda_N^{a_*}$. 
There exists $\eta_*>0$ such that, for every sufficiently small $\tau>0$, $\delta>0$, and for any smooth functions
$\eta$, $I$ and $a$ such that, for all $t\in[0,\tau]$
$$\eta(t)\geq\eta_*~,~~ I(t)\geq 0~,~~ a_*-\delta =a(0) \leq a(t) \leq a(\tau)=a^*+\delta~\text{ and }~|a'(t)|\leq \kappa,$$
we have
$$\forall k\leq N~~,~~~\|\G_{0}^{\tau}\phi_k^{\infty,a_*-\delta}-\phi_{\sigma(k)}^{\infty,a_*+\delta}\|_{L^2}\leq \varepsilon,$$
where $\G_s^t$ the flow of the Schr\"odinger equation \eqref{eq}, $\sigma=\sigma_{a_*-\delta}^{a_*+\delta}$ is the permutation described in Section \ref{intro},
and $\phi_k^{\infty,a}$ are the eigenmodes of $H^{\infty, a}$.
\end{prop}

{\noindent \emph{\textbf{Proof}:}}
Note that the eigenmodes $\phi_k^{\infty,a(t)}$ given by \eqref{def_varphi_l},\eqref{def_varphi_r}
are not of class $H^2$ in the whole interval $(0,1)$ and it does not vanish everywhere in 
$\big[a(t)-\frac{1}{\eta_*},a(t)+\frac{1}{\eta_*}\big]$, so we cannot apply Lemma \ref{lemma-compare}
to these functions directly. We therefore modify them near the support interval of $V$. To fix the notations, 
choose some $k\leq N$ and assume that $\phi_k^{\infty,a_*-\delta}$ is given by \eqref{def_varphi_l} for some $p\in\NN^*$ (the estimates for the case where $\phi_k^{\infty,a_*-\delta}$ is given by \eqref{def_varphi_r} are similar). Let $\chi$ be a smooth truncation equal to $0$ in $[-1,+\infty)$ and to $1$ in $(-\infty,-2]$. We denote $\chi_\alpha(\cdot)=\chi(\cdot/\alpha)$; this function vanishes in $[-\alpha,+\infty)$. Denote
\begin{equation}\label{def_psi}
\psi_k^{a(t)}(x):= \sqrt{\frac{2}{a(t)}}\sin\Big(\frac{p\pi}{a(t)} x \Big)\chi_\alpha(x-a(t)).
\end{equation}
Now, $\psi_k^{a(t)}$ is a smooth function vanishing in 
$\big[a(t)-\frac{1}{\eta_*},a(t)+\frac{1}{\eta_*}\big]$ when $\alpha\geq\frac{1}{\eta_*}$. Moreover, $\psi_k^{a(t)}$ is close, uniformly for all $t$, to $\phi_k^{\infty,a(t)}=\sqrt{\frac{2}{a(t)}}\sin\big(\frac{p\pi}{a(t)} x \big)\Un_{x\in[0,a(t)]}$ in $L^2((0,1),\CC)$ when $\alpha$ is small. More precisely, 
\begin{align}\label{eq_smoothtruncation}
\|\phi_k^{\infty,a }-\psi_k^{a}\|_{L^2}&
\leq\sqrt{\frac{2}{a}\int_{a-2\alpha}^a \sin^2\big(\frac{p\pi}{a}x\big)\dd x }
\leq \frac{4\pi}{\sqrt{3}} p \left(\frac{\alpha}{a}\right)^{\frac{3}{2}} = \Oc (N\alpha^{3/2}). 
\end{align}

To apply Lemma \ref{lemma-compare}, we need to estimate the derivatives of $\psi_k^{a(t)}$ when 
$\alpha$ is small. Notice that the first and second derivatives of $\chi_\alpha$ are 
of order $\frac{1}{\alpha}$ and, respectively, $\frac{1}{\alpha^2}$, and are supported 
in $[a(t)-2\alpha,a(t)-\alpha]$. Thus, we have 
\begin{align*}
\int_0^1\Bigg|\sin\Big(\frac{p\pi}{a(t)} x \Big)\ddd^2_{xx}\chi_\alpha(x-a(t))\Bigg|^2\dd x  &\leq \frac{p^2\pi^2\|\ddd^2_{xx}\chi\|_\infty^2}{\alpha^4a(t)^2}\int_{a(t)-\alpha}^{a(t)-2\alpha}(a(t)-x)^2\dd x\\ 
&= \Oc(N^2 \alpha^{-1})
\end{align*}
and, similarly, we obtain that 
$\int_0^1|\cos(\frac{p\pi}{a(t)} x )\ddd_{x}\chi_\alpha(x-a(t))|^2\!\dd x = 
\Oc(\alpha^{-1})$ and that $\int_0^1|\sin(\frac{p\pi}{a(t)} x )\chi_\alpha(x-a(t))|^2\!\dd x = 
\Oc(N^2\alpha^{3})$. This gives
\begin{align}\label{estimations}
\left\|\partial_{xx}^2\psi_k^{a(t)}\right \|_{L^2} \; \leq & \;\;
\sqrt{\frac{2}{a(t)}}\; \left[\ \left\|\sin\Big(\frac{p\pi}{a(t)}\ \cdot \Big) \
\ddd^2_{xx}\chi_\alpha(\ \cdot \ - a(t)) \right\|_{L^2} \right. \nonumber \\ 
& \left. +\
2\frac{p\pi}{a(t)} \left\|\cos\Big(\frac{p\pi}{a(t)}\ \cdot \!\Big)\
\ddd_{x}\chi_\alpha(\ \cdot \ - a(t)) \right\|_{L^2} \right. \nonumber \\ 
& \left. +\
\frac{p^2\pi^2}{a(t)^2} \left\|\sin\Big(\frac{p\pi}{a(t)}\ \cdot \!\Big) \chi_\alpha(\ \cdot \ - a(t)) \right\|_{L^2}\right]
\nonumber\\
& =\Oc(N\alpha^{-1/2}+N^3\alpha^{3/2}).
\end{align}
By similar computations, 
\begin{equation}\label{estimationst}
\left \|\partial_t \psi_k^{a(t)} \right\|_{L^2}  = \Oc \left(N\sup_{t\in[0,t]} |a'(t)| \right) = \Oc(N\kappa).
\end{equation}
Notice that the $L^2-$norm of $\psi^a_k$ is controlled by $\Oc(1+N\alpha^{3/2})$ due to \eqref{eq_smoothtruncation}. Now, by Lemma \ref{lemma-compare}, if
$u(\tau)=\Gamma^\tau_{0}\phi_k^{\infty,a_*-\delta}$, then, by \eqref{estimations} and \eqref{estimationst},
$$\|u(\tau)-\psi_k^{a(\tau)}\|_{L^2}^2 \leq \|\phi_k^{\infty,a_*-\delta}-\psi^{a_*-\delta}_{k}\|_{L^2}^2 +
N \tau \; \Oc(\alpha^{-1/2} + N^2\alpha^{3/2} + (2+N\alpha^{3/2})\kappa).$$
Using \eqref{eq_smoothtruncation}, we can then estimate
\begin{align}\label{eq_overall}
\!\!\!\!\!\!\!\|u(\tau)-&\phi_{\sigma(k)}^{\infty,a(\tau)}\|_{L^2} \leq  \|u(\tau)-\psi^{a(\tau)}_k\|_{L^2}+\|\psi^{a(\tau)}_k-\phi_{\sigma(k)}^{\infty,a(\tau)}\|_{L^2} \nonumber \\
&\leq\|\phi_k^{\infty,a_*-\delta}-\psi^{a_*-\delta}_k\|_{L^2}
+
\sqrt{N\tau} \; \Oc\left(\sqrt{\alpha^{-1/2} + N^2\alpha^{3/2} + (2+N\alpha^{3/2})\kappa}\;\right) \nonumber \\ &\qquad\qquad~~~~+\|\psi^{a_*+\delta}_k-\phi_{\sigma(k)}^{\infty,a_*+\delta}\|_{L^2} \nonumber \\
&= \Oc (N\alpha^{3/2}) + \sqrt{N\tau} \; \Oc\left(\sqrt{\alpha^{-1/2} + N^2\alpha^{3/2} + (2+N\alpha^{3/2})\kappa}\;\right).
\end{align}
Choose $\alpha \ll \left(\varepsilon/N\right)^{2/3}$. Set $\eta_*>0$ sufficiently large so that 
$\alpha\geq \frac{1}{\eta^*}$. Then, the above computations are valid for any choice of the function $\eta(t)$ bounded by $\eta_*$ from below for all $t\in [0,\tau]$. Therefore, if we take
$$\tau \ll \frac{\varepsilon^2}{N (\alpha^{-1/2} + N^2\alpha^{3/2} + (2+N\alpha^{3/2})\kappa)},$$
the estimate \eqref{eq_overall} yields the claim. Note that the only restriction on $\delta$ (the
range of the displacement of the potential wall) is given by 
$$\delta \leq \sup_{t\in[0,t]} |a'(t)| \tau \leq \kappa \tau,$$
so the result holds true for all sufficiently small $\delta$.
$\hfill\square$

\vspace{5mm}

{\noindent \bf Remark:}
From \eqref{eq_overall}, we may estimate the behavior of the parameters in the statement of Proposition \ref{prop-compare-eigenfunctions1}:
\begin{itemize}
\item For $N\in\NN^*$ and $\kappa>0$ fixed, when the error $\varepsilon$ is small, the sharpness of the potential $\eta^*$ is at least $\Oc(\epsilon^{-\frac{2}{3}})$. Hence, $\tau$ and $\delta$ have to be at most $\Oc(\epsilon^\frac{7}{3})$.
\item For $\epsilon>0$ and $\kappa>0$ fixed, when we consider a large number $N$ of frequencies, then $\eta^*$ is at least $\Oc(N^{\frac{2}{3}})$. Thus, $\tau$ and $\delta$ have to be at most $\Oc(N^{-2})$. 
\item Fix $\epsilon>0$ and $N\in\NN^*$. On the one hand, if we choose a slow motion and $\kappa$ is small, then the values of $\eta_*$ and $\tau$ are not affected by $\kappa$, while the distance $\delta$ is of order $\Oc(\tau \kappa)=\Oc(\kappa)$. On the other hand, when $\kappa$ is large, the time $\tau$ is of order $\Oc(1/\kappa)$ (again $\eta_*$ and $\delta$ do not depend on $\kappa$).
\end{itemize}


\section{Proof of Theorem \ref{th}}\label{section-proof-of-theorem} 
\begin{figure}[ht]
\resizebox{0.95\textwidth}{!}{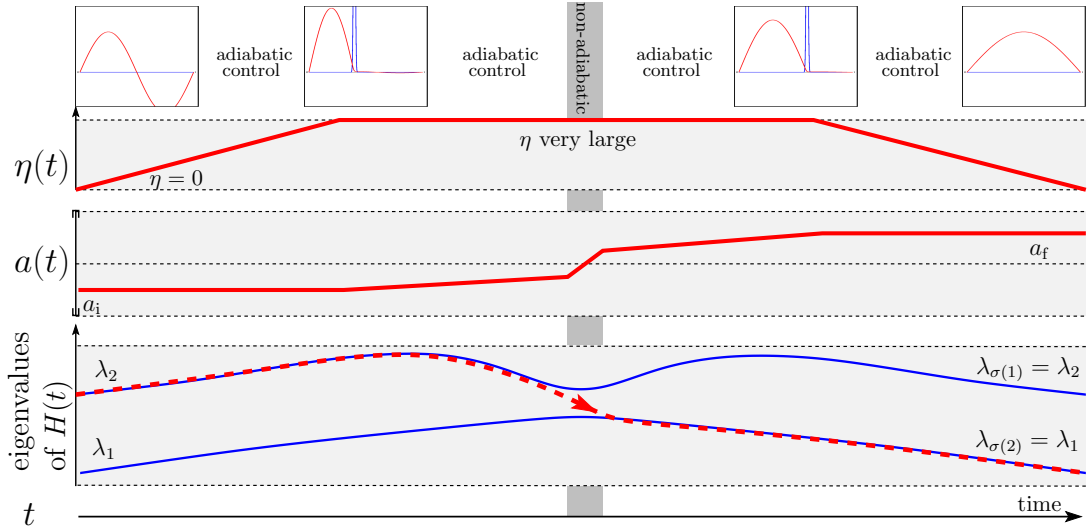}
\caption{\it The figure represents the control strategy behind Theorem \ref{th}, steering $\sin(2\pi x)$ in $\sin(\pi x)$. In the first step, the potential adiabatically increases at the fixed location $\ai$. During the second step, the location of the potential wall $a(t)$ is moved from $\ai$ to $\af$. Close to $a_*=1/2$, we have a short-time transition from $\lambda_2$ to $\lambda_1$. In the final step, the potential adiabatically vanishes.}\label{fig_preuve}
\end{figure}

In order to obtain the permutation of eigenmodes described in Theorem \ref{th}, we use a control path defined as follows (see also Figure \ref{fig_th}). We fix some very large $\eta_*$ and, first, start to adiabatically increase $I$ from zero to a very large value $I_*$, thus obtaining a very thin and high potential wall located at $\ai$. The evolution of eigenstates is then described by the classical adiabatic result stated in Proposition \ref{prop-adiab}. Second, we move the wall to the location $\af$. In this step, the potential $V(x,t)=I_*\rho^{\eta_*}(x-a(t))$ has a large fixed amplitude and a small fixed width $(\eta_*)^{-1}$,
while the position $a(t)$ of the wall is moving. The movement of the wall alternates adiabatic motion (by Proposition \ref{prop-adiab}) and short-time level crossings (by 
Proposition \ref{prop-compare-eigenfunctions1}). Finally, we adiabatically decrease $I$ to zero, up to the
extinction of the potential. For further details, we refer to Figure \ref{fig_preuve}.

\vspace{3mm}

\noindent {\bf The crossings of eigenvalues and the corresponding permutations.} To simplify the notations, we assume that $\ai<\af$ and we consider any $a(t)$ monotonically going from $\ai$ to $\af$, such that
$|a'(t)|\leq \kappa$ for all $t$ (where $\kappa$ is the bound from the statement of Theorem \ref{th}). 

The operator $H^{\infty, a}$ (the Laplacian with the Dirichlet boundary conditions on the split interval
$(0,a)\cup (a,1)$, as introduced in Section \ref{section_spectre}), corresponds to ideal motion described by the formal equation \eqref{eq_dmitry}. We recall that its eigenvalues $(\lambda^{\infty,a}_k)_{k\in\NN^*}$ are grouped into two families $(\mu^l_p(a))_{p\in\NN^*}$ and $(\mu^r_q(a))_{q\in\NN^*}$ given by \eqref{def_mu}. During the motion of $a$ from $\ai$ to $\af$, it can happen that some spectral curves belonging to $(\mu_p^l(a))_{p\in\NN^*}$ cross some others in $(\mu_q^r(a))_{q\in\NN^*}$ and vice versa. For $N\in\NN^*$ (where $N$ is given in the statement of Theorem \ref{th}), 
there exist $P,Q\in\NN$ so that 
$$(\lambda_k^{\ai})_{k\leq N}\ =\ \big(\mu_p^l(\ai)\big)_{p\leq P}\ \cup\ \big(\mu_q^r(\ai)\big)_{q\leq Q}~.$$
When $a$ moves from $\ai$ to $\af$, the spectral curves respectively belonging to $\big(\mu_p^l({a})\big)_{p\leq P}$ and $\big(\mu_q^r({a})\big)_{q\leq Q}$ can cross a finite number of other curves at a finite number of locations of $a$. We call these locations $(a_{*,j})_{j\leq J}$ with 
$J\in\NN^*$ and
$$\ai~<~a_{*,1}~<~\ldots~<~a_{*,J}~<~\af~$$
(as $\ai$ and $\af$ are irrational, no crossing can happen at the end points $a=\ai$ or $a=\af$). We denote by $M\in\NN^*$ the smallest number such that any of the curves $(\mu_p^l(a))_{p\leq P}$ and $(\mu_q^r(a))_{q\leq Q}$ cross only eigenvalues in $(\lambda^{\infty,a}_k)_{k\leq M}$ when $a$ moves from $\ai$ to $\af$. As in Theorem \ref{th}, for any $a$ and $b$ irrational in $[\ai,\af]$, we define the permutation $\sigma_a^b$  as follows:
\begin{itemize}
\item if $k\in\NN^*$ is such that $\lambda^{\infty,a}_k$ is equal to $\mu_p^l(a)=\frac{p^2\pi^2}{a^2}$ for some $p$, then $\sigma_a^b(k)$ is the index such that $\mu_p^l(b)=\frac{p^2\pi^2}{b^2}=\lambda_{\sigma_{a}^{b}(k)}^b$;
\item if $k\in\NN^*$ is such that $\lambda^{\infty,a}_k$ is equal to 
$\mu_q^r(a)=\frac{q^2\pi^2}{(1-a)^2}$ for some $q$, then $\sigma_a^b(k)$ is the index such that $\mu_q^r(b)=\frac{q^2\pi^2}{(1-b)^2}=\lambda_{\sigma_{a}^{b}(k)}^b$. 
\end{itemize}
Let $\delta_*>0$ be small enough, so that $\ai<a_{*,1}-\delta_*$, $a_{*,j}+\delta_*<a_{*,j+1}-\delta_*$ for $j=1,\ldots,J-1$ and $a_{*,J}+\delta_*<\af$. As $(a_{*,j})_{j\leq J}$ correspond to the crossings point of the $N$ first eigenvalues, we have 
$$\forall k\leq N~,~~\sigma_\ai^\af(k)~=~\sigma_{a_{*,J}-\delta_*}^{a_{*,J}+\delta_*}\,\circ\,\ldots\,\circ\,\sigma_{a_{*,1}-\delta_*}^{a_{*,1}+\delta_*}\,(k)~.$$

\vspace{3mm}

\noindent {\bf Fixing the parameters.} Let $\varepsilon'<\varepsilon/(4J+3)$, where $\varepsilon>0$ is the small error introduced in Theorem \ref{th}. 
We apply Proposition \ref{prop-compare-eigenfunctions1} to each $a_*=a_{*,j}$ where crossings of eigenvalues occurs. We obtain some large $\eta_*$ and we can choose, for every $j\leq J$, a distance $\delta_j<\delta_*$, a time $\tau_j$ and a path $a_j\in\Cc^\infty([0,\tau_j],[a_{*,j}-\delta_j,a_{*,j}+\delta_j])$ satisfying 
\begin{equation}\label{conditions_on_the_jump}\|a'\|_{L^\infty([0,\tau],\RR)}\leq \kappa,\ \ \ \ \ \ a_j^{(m)}(0)=a_j^{(m)}(\tau_j)=0,\ \ \ \ \ \ \forall m\in\NN^*\end{equation} 
such that 
\begin{equation}\label{eq_demo_1}
\forall \eta\geq \eta_*~,~~
\forall k\leq N~,~~\|\Gamma^\tau_0 \phi_k^{\infty,a_{*,j}-\delta_j} - \phi_{\sigma_j(k)}^{\infty,a_{*,j}+\delta_j}\|_{L^2}\leq \varepsilon'~.
\end{equation}
Next, we apply Theorem \ref{th-cv-spectre} for each $a=a_{*,j}\pm\delta_j$ and for $k\leq M$. We obtain that, maybe with a larger $\eta_*$, the estimate \eqref{eq_demo_1} is also guaranteed for
$\eta = \eta_*$ and $I=I_*$  (for some sufficiently large $I_*$):
\begin{equation}\label{eq_demo_2}
\forall k\leq M,~\exists \alpha_k^{I_*,\eta_*,j\pm}\in\UU,~~\|\phi_k^{I_*\eta_*,a_j\pm\delta_j}-\alpha_k^{I_*,\eta_*,j\pm}\phi_k^{\infty,a_j\pm\delta_j}\|_{L^2} \leq \varepsilon'
\end{equation}
where $\UU=\{z\in\CC~,~|z|=1\}$ contains all the possible phase-shifts.

\vspace{3mm}

\noindent {\bf The initial ``vertical'' adiabatic motion.} We start to construct the potential $V(x,t)=I(t)\rho^{\eta(t)}(x-a(t))$ as follows. First, we fix $a(t)$ to be constantly equal to $\ai$, $\eta(t)=\eta_*$, and we choose $I(t)$ to be a smooth function going from $0$ to $I_*$, with 
the derivative $I'(t)$ compactly supported in $(0,T_1)$ and satisfying $|I'(t)|<\kappa$ for some sufficiently large $T_1>0$. Then, we apply 
Proposition \ref{prop-adiab} and obtain an adiabatic evolution lasting up to the time $t=T_1$
such that 
\begin{equation}\label{eq_demo_3}
\forall k\leq N~,~~\exists \alpha_k\in\UU~,~~\|\Gamma_0^{T_1} \sqrt{2}\sin(k\pi\cdot)-\alpha_k \phi_k^{I_*,\eta_*,\ai}\|_{L^2} \leq \varepsilon'~.
\end{equation}

\vspace{3mm}

\noindent {\bf The initial ``horizontal'' adiabatic motion.} Now, we fix $I(t)$ to be a constant function equal to $I_*$ and we choose $a(t)$ to be a smooth monotone function going from $\ai$ to $a_{*,1}-\delta_1$, with the derivative $a'$ compactly supported in $(0,T_2)$ for a sufficiently large $T_2>0$. We apply again Proposition \ref{prop-adiab} and we obtain an adiabatic motion lasting the time $T_2$ such that  
\begin{equation}\label{eq_demo_4}
\forall k\leq N~,~~\exists \alpha_k\in\UU~,~~\|\Gamma_0^{T_2} \phi_k^{I_*,\eta_*,\ai}-\alpha_k \phi_k^{I_*,\eta_*,a_{*,1}-\delta_1}\|_{L^2} \leq \varepsilon'~.
\end{equation}
As above, we choose $T_2$ sufficiently large so that $\|a'\|_{L^\infty((0,T_2),\RR)}\leq \kappa.$ Now, we concatenate the obtained path with the previous one. We notice that the concatenation is still smooth since the derivatives vanish near the joint-point. Moreover, $\Gamma_s^t$ is unitary, so the $L^2-$errors propagate without changing. Gathering \eqref{eq_demo_3} and \eqref{eq_demo_4}, we obtain the paths $a(t)$, $I(t)$, and $\eta(t)$, both smooth in $[0,T_1+T_2]$, satisfying   
\begin{equation}\label{eq_demo_5}
\forall k\leq N~,~~\exists \alpha_k\in\UU~,~~\|\Gamma_0^{T_1+T_2} \sqrt{2}\sin(k\pi\cdot)-\alpha_k \phi_k^{I_*,\eta_*,a_{*,1}-\delta_1}\|_{L^2} \leq 2\varepsilon'~.
\end{equation}
In addition, we still have that $\|I'\|_{L^\infty((0,T_1+T_2),\RR)}\leq \kappa$, 
$\|\eta'\|_{L^\infty((0,T_1+T_2),\RR)}\leq \kappa$, and $\|a'\|_{L^\infty((0,T_1+T_2),\RR)}\leq \kappa$.

\vspace{3mm}

\noindent {\bf The first short-time crossing.} In the previous step, we stopped at $a_{*,1}-\delta_1$, just before the first crossing point of the ideal eigenvalues $\lambda^{\infty,a}_k$. At this point, we can no longer follow the perfectly adiabatic motion (as explained in Section \ref{section_fast_translation}), because the real eigenvalues $\lambda_k^{I,\eta,a}$ is simple for every $a\in (0,1)$ (see Figure \ref{fig_crossing}). Thus, we use a quasi-adiabatic path $a_1(t)$ lasting a time $T_3$ to go from $a_{*,1}-\delta_1$ to $a_{*,1}+\delta_1$ satisfying \eqref{conditions_on_the_jump} and \eqref{eq_demo_1}. We concatenate this path to the previous ones. Once again, the evolution operator is unitary and the error term from \eqref{eq_demo_5} is transmitted without amplification. Due to \eqref{conditions_on_the_jump}, \eqref{eq_demo_1} and \eqref{eq_demo_2}, the concatenated control satisfies
\begin{equation*}
\forall k\leq N~,~~\exists \alpha_k\in\UU~,~~\|\Gamma_0^{T_1+T_2+T_3} \sqrt{2}\sin(k\pi\cdot)-\alpha_k \phi_{\sigma_1(k)}^{I_*,\eta_*,a_{*,1}+\delta_1}\|_{L^2} \leq 5\varepsilon'~.
\end{equation*}
Again we keep the controls $\|I'\|_{L^\infty((0,T_1+T_2+T_3),\RR)}\leq \kappa$, 
$\|\eta'\|_{L^\infty((0,T_1+T_2+T_3),\RR)}\leq \kappa$,
and $\|a'\|_{L^\infty((0,T_1+T_2+T_3),\RR)}\leq \kappa$.

\vspace{3mm}

\noindent {\bf Iteration of the process.} We repeat the strategy presented above for every crossing point. In particular, we proceed as follows for every $1\leq j\leq J$. First, we move adiabatically between 
$a_{*,j}+\delta_j$ and $a_{*,j+1}-\delta_{j+1}$, which add an error $\varepsilon'$, like in \eqref{eq_demo_4}. This motion preserves the ordering of the eigenvalues. The duration of the motions is chosen large enough in order to control the speed of the wall's motion with the parameter $\kappa$. Second, we have the quasi-adiabatic dynamics satisfying \eqref{conditions_on_the_jump}, \eqref{eq_demo_1} and
\eqref{eq_demo_2} while moving from $a_{*,j+1}-\delta_{j+1}$ to $a_{*,{j+1}}+\delta_{j+1}$. This motion adds an error $3\varepsilon'$ and applies the permutation $\sigma_{j+1}$ to the eigenmodes. Finally, we move adiabatically between 
$a_{*,J}+\delta_j$ and $\af$ as in the first ``horizontal'' adiabatic motion by controlling again the speed of the wall's motion. We obtain that, for a suitable time $T$, there exist smooth controls $I(t)$, $\eta(t)\equiv\eta_*$, and $a(t)$ such that 
$$\|\eta'\|_{L^\infty((0,T),\RR)}\leq \kappa, \qquad \|\eta'\|_{L^\infty((0,T),\RR)}\leq \kappa,\ \ \ \ \ \ \ \ \ \ \ \ \|a'\|_{L^\infty((0,T),\RR)}\leq \kappa,$$
\begin{equation*}\begin{split}
\forall k\leq N~,~~\exists \alpha_k\in\UU~,~~\|\Gamma_0^{T} \sqrt{2}\sin(k\pi\cdot)-\alpha_k \phi_{\sigma(k)}^{I_*,\eta_*,a_{f}}\|_{L^2} \leq (4J+2)\varepsilon'~
\end{split}\end{equation*}
where $\sigma=\sigma_{J}\circ\ldots\circ\sigma_1$ is exactly the permutation $\sigma_\ai^\af$ described in  Theorem \ref{th}. 

\vspace{3mm}

\noindent {\bf Final ``vertical'' adiabatic motion.} With fixed $a(t)\equiv \af$, and $\eta(t)=\eta_*$, we adiabatically decrease $I$ to $0$ with the speed slower than $\kappa$. The whole concatenated control path gives us Theorem \ref{th}, as
\begin{equation*}
\forall k\leq N~,~~\exists \alpha_k\in\UU~,~~\|\Gamma_0^{T} \sqrt{2}\sin(k\pi\cdot)-\alpha_k \sqrt{2}\sin(\sigma(k)\pi\cdot)\|_{L^2} \leq (4J+3)\varepsilon' <\varepsilon~.
\end{equation*}


\section{Further discussions}\label{further_discussions}
\subsection{The techniques behind Theorem \ref{th} in practice}\label{applications}

An interesting peculiarity of the strategy employed in Theorem \ref{th} is that every aspect of it can be made explicit. In particular, by retracing each part of the proof,
one can provide explicit controls $a$, $I$, and $\eta$, and also specify the time $T$ in terms of
the data $a_i,$ $a_f$, $N$, $\varepsilon$, and $\kappa$. To this purpose, the first step consists in evaluating the parameters involved in the statement of Proposition \ref{prop-compare-eigenfunctions1}. Second, one has to specify the rate of convergence given by Theorem \ref{th-cv-spectre} and state 
Proposition \ref{prop-adiab} with explicitly given time interval and control path. Then all these steps should be implemented as it is done in the proof of Theorem \ref{th}. The corresponding result would be interesting both from a theoretical point of view and for some practical implementations. Indeed, numerical simulations could be done to reproduce the results of this work and characterize the different phenomena presented. For instance, one could be interested in studying the transition of the eigenmodes from the right side of the potential wall $V$ to the left one (or vice versa) described by Figure \ref{fig_crossing} by quantifying the tunnel effect in terms of the parameters $\eta$ and $I$.
 
\vspace{3mm}
 
Another interesting aspect is that we have plenty of freedom in the choice of the control path, which can be modified according to other preferences. First, all the adiabatic paths of the strategy can be varied freely, except concerning the starting point and the final point and the fact that they should be slow enough to apply Proposition \ref{prop-adiab}. Second, the quasi-adiabatic dynamics defined by Proposition \ref{prop-compare-eigenfunctions1}, can be modified too. Indeed, this Proposition allows to substitute the wall
movement constructed in the proof with other types of motion. For instance, we can let the potential located at $a_*-\delta$ decrease non-adiabatically up to extinction and, afterwards, grow non-adiabatically the potential wall at $a_*+\delta$. Both motions have to be fast enough; they can be obtained by varying the value of the intensity $I$, which does not interfere with the validity of Proposition \ref{prop-compare-eigenfunctions1}. Another type of dynamics, which could be even more interesting for practical implementations, is achieved by an instantaneous displacement of the potential wall from the $a=a_*-\delta$ to $a=a_*+\delta$. Even though such kind of dynamics would provide a discontinuous control, it would lead to the same outcome of Proposition \ref{prop-compare-eigenfunctions1} and it would be simpler for numerical implementations. 

\subsection{Exponential energy growth}\label{energy_growth}
It was noted in \cite{Dmitry} that the repetition of the permutation described in Theorem \ref{th} should, typically, lead to an exponential growth in energy. Namely, if we start with an eigenstate $\psi_k$ with the 
eigenvalue $\lambda_k$, then after applying $n$ times the permutation $\sigma_{\ai}^{\af}$ we
will find the system at the eigenstate $\psi_{k_n}$ with $k_n = \left(\sigma_{\ai}^{\af}\right)^n (k)$.
The claim is that the energy $\lambda_{k_n}$ behaves as
$$\lambda_{k_n} \sim  e^{n r} \lambda_k$$
for some $r(\ai,\af)>0$, for a typical initial condition $k_0=k$ and typical values of $\ai$ and $\af$.

Proving this (after a proper rigorous reformulation) could be a non-trivial task. Still, it is not hard to build
explicit examples of the exponential growth, see \cite{Dmitry}. The argument behind the general
claim is as follows. Note that for any permutation $\sigma$ of the set of natural numbers its trajectory $k_n=\sigma^n (k_0)$ is either looped or tends to infinity both at forward and backward iterations.
The exponential growth claim means, therefore, that for the permutations $\sigma_{\ai}^{\af}$ we consider here the most of trajectories are not looped and the averaged value of $\ln \frac{k_{n+1}}{k_n}$ is strictly positive typically (as $\lambda_k\sim k^2$, the
exponential growth of the energy is equivalent to the exponential growth of the eigenstate number $k$). 

In order to estimate the average value of the increment in $\ln k$, we recall that at each irrational $a$
the eigenstates of the Laplacian in the split interval $(0,a)\cup (a,1)$ with the Dirichlet boundary conditions
are divided into two groups, the left eigenstates are supported in $(0,a)$ and the right ones are supported
in $(a,1)$. Thus, given irrational $a$, we can order the eigenstates by their energy $\lambda_k$ and introduce the indicator sequence: $\xi^a (k)=1$ if the eigenstate $\psi_k$ is left, and $\xi^a(k)=-1$
if $\psi_k$ is right. The two sequences $\xi^{\ai}$ and $\xi^{\af}$ completely determine the permutation $\sigma_{\ai}^{\af}$. Indeed, if at some $a$ the state $\psi_k$ is left and acquires the number $m$
when we order the left states by the increase of energy, then there are exactly $m$ left and 
$(k-m)$ right states with the energies not exceeding $\lambda_k$, so
$$S^a(k) := \xi^a(1) + \dots + \xi^a (k) = 2m - k;$$
if $\psi_k$ is a right state with the number $n$ in its group, then there are exactly $n$
right and $(k-n)$ left states with the energies not exceeding $\lambda_k$, so 
$$S^a(k) = k − 2n.$$
These two formulas can be rewritten as
$$m \mbox{  or  } n =\frac{1}{2} (k + \xi^a(k) S^a(k)).$$

Thus, after the splitting the interval $(0,1)$ at $a=\ai$ the state $\psi_k$
becomes left if $\xi^{\ai}(k)=1$ or right if $\xi^{\ai}(k)=-1$. When $a$ changes from
$\ai$ to $\af$ left states remain left and right states remain right, and the corresponding number 
$m$ or $n$ stays constant, implying that the number $k$ changes to $\bar k= \sigma_{\ai}^{\af}$ 
according to the rule
\begin{equation}\label{rulek}
k + \xi^{\ai}(k) S^{\ai}(k) = \bar k + \xi^{\ai}(k) S^{\af}(\bar k).
\end{equation}
One needs a proper analysis of dynamics of $k$ defined by this formula, but we just make a heuristic assumption that $\xi^{ai}$ and  $\xi^{af}$ are sequences of independent, identically distributed
random variables. Let $\beta \in(0,1)$ be the probability of $\xi^{\ai} = 1$ and $\gamma\in (0,1)$
be the probability of $\sigma^{\af}=1$. Then, in the limit of large $k$,
we may substitute $S^{\ai}(k) \sim (2\beta-1)k$ and $S^{\af}(\bar k) \sim (2\gamma-1)\bar k$ in (\ref{rulek}), which gives
$$\bar k =\left\{\begin{array}{ll} \frac{\beta}{\gamma} k +o(k) & \mbox{ with probability  } \beta,\\
\frac{1-\beta}{1-\gamma} k +o(k) & \mbox{ with probability  } 1  - \beta. \end{array}\right.$$
It follows that
$$r := \mathbb{E}(\ln \bar k - \ln k) =\beta \ln\frac{\beta}{\gamma} + 
(1-\beta) \ln\frac{1-\beta}{1-\gamma}$$
in the limit of large $k$. This quantity is strictly positive if $\beta\neq \gamma$, giving the claimed exponential growth. It would be interesting to do a more rigorous analysis of the dynamics generated
by permutations of eigenvalues due to this or other cyclic quasi-adiabatic processes.

The process described by equation \eqref{eq}, which we study in this paper, follows the ``ideal'' permutation
$\sigma_{\ai}^{\af}$ only approximately, and only until the eigenstate number $k$ is smaller than a certain fixed number $N$. One can check through the proofs that one can indeed repeat this process
until the corresponding eigenvalue remains significantly lower than $I_*$, the maximal intensity of the potential $V(x,t)= I(t) \rho^{\eta(t)}(x-a(t))$. Thus, until this moment, we can expect the exponential growth of the energy at the repeated application of the control cycle described in Theorem \ref{th}.
This lets us to estimate as $\Oc(\ln I_*)$ the number of our control cycles sufficient to transform a low energy state to the state with energy of order $I_*$. As the speed with which we change $I$ is bounded from above, and we have $I=0$ at the beginning of each cycle, the duration of one cycle is $\Oc(I_*)$.
Thus, we conclude that we can transform a low energy eigenstate to the state of energy of order $I_*$
in a time of order $I_*\ln I_*$.

In fact, no unbounded exponential growth of the energy is possible in our setting, or for any 
periodic in time, smooth potential $V(x,t)$ in a bounded domain. Indeed, assume that $\psi(t)$ is a solution of the Schr\"odinger equation with $\psi(0)$ concentrated close to the $k-$th eigenmode. Until our permutation process stay valid, after $n$ cycles, $\psi(nT)$ is concentrated close to the $\sigma^n(k)-$th eigenmode. This implies in particular that $\|\psi(nT)\|_{H^1}$ grows like $\Oc(\sigma^n(k))$. However, it is proved in \cite{Bourgain} (see also \cite{Delort}) that the $H^s-$norm of the solution of a periodic Schr\"odinger equation has at most a polynomial growth rate for any $s>0$. Even if the ideal permutation $\sigma$ yields an exponential growth of the eigenmodes, we must move away from this aimed trajectory after a finite time. The realization of the exponential growth needs to use a more singular evolution as the splitting process introduced in \cite{Dmitry}.

\subsection{Removing the phase shifts in Theorem \ref{th}}\label{section_phase}

We notice that permutation of the eigenmodes described in Theorem \ref{th} is obtained up to phase shifts. In order to explain how to remove them, we present the following example. Let $\alpha_1\phi_1^{I,\eta,a}$ and $\alpha_2\phi_2^{I,\eta,a}$ with $\alpha_1,\alpha_2\in \UU$ be two eigenstates of $-\partial_{xx}^2+I\rho^{\eta}(x-a)$ 
with $I>0$, $\eta>0$ and $a\in (0,1)$ fixed. We assume that the two corresponding eigenvalues are rationally independent. If we apply the dynamics generated by such Hamiltonian and we wait a time $T$, then the two states respectively become $e^{-iT\lambda_1}\alpha_1 \phi_1^{I,\eta,a}$ and 
$e^{-iT\lambda_2}\alpha_2 \phi_2^{I,\eta,a}$. Fixing an error $\varepsilon>0$, the rational independence ensures that there is a time $T$ such that 
\begin{equation}\label{phases} e^{-iT\lambda_1}\alpha_1 \phi_1^{I,\eta,a}=
\phi_1^{I,\eta,a}+\Oc(\varepsilon),\ \ \ \ \ \  \ 
e^{-iT\lambda_2}\alpha_2 \phi_2^{I,\eta,a}= \phi_2^{I,\eta,a}+\Oc(\varepsilon).\end{equation}
The same is valid when we consider $N$ eigenmodes and we assume that the corresponding eigenvalues are rationally independent. As explained in \cite[Proposition\ 3.2]{Marion}, the property of being rationally independent for the eigenvalues of an operator $-\partial_{xx}^2+V$ with $V\in L^\infty((0,1),\RR)$ is generically satisfied with respect to $V$. 

Thus, if we want to obtain Theorem \ref{th} with an error $\varepsilon>0$ and without phases, consider a control path constructed in the proof of Theorem \ref{th} for an error $\varepsilon/3$. If we slightly perturb $\rho$, the result is unchanged, up to an additional error $\varepsilon/3$. Thus, we can consider $\rho$ such that the eigenvalues of $-\partial_{xx}^2+V$ are irrationally independent for $V=I\rho^\eta(\cdot-a)$ with $I$, $a$ and $\eta$ being some values of the parameters reached during the control path. At this moment of the control path, it is sufficient to pause and wait long enough in order to eliminate the phases
up to a small error $\varepsilon/3$, as done in \eqref{phases}. 

Notice that it is, in fact, possible to adjust the phase shifts without assumptions on the spectrum of $-\partial_{xx}^2+V$ if $V$ is composed by several walls as in Theorem \ref{th2} and Theorem \ref{th3}. The idea is explained in Section \ref{section_approx} below.

\subsection{Realizing an arbitrary finite permutation}\label{any_permutation}

It is clear that the arguments of the proof of Theorem \ref{th} can be generalized to the case where the potential is given by 
\begin{equation}\label{eq_potentiel_3}
V(t,x)~=~\sum_{j=1}^J  I_j(t) \rho^{\eta_j(t)}(x-a_j(t))  
\end{equation}
with several moving walls which ``almost split'' the interval. By changing the walls positions, we can control the state of the system by the quasi-adiabatic motion. Theorem \ref{th} allows us to realize a specific permutation $\sigma_\ai^\af$ to the eigenmodes. It is not clear which permutations are reachable by combining several permutations of this specific type. However, with $V$ given by \eqref{eq_potentiel_3}, it is easy to check that we can reach any permutation of any given number of the low modes.
\begin{prop}\label{th2}
Let $\sigma:\NN^*\rightarrow \NN^*$ be any permutation. For all $N\in\NN^*$ and $\varepsilon>0$, 
there exist $J\in\NN$, $T>0$, and smooth functions 
$I_j\in\Cc^\infty([0,T],\RR^+)$,
$\eta_j\in\Cc^\infty([0,T],\RR^+)$, and $a_j\in\Cc^\infty ([0,T],(0,1))$, $j=1\ldots J$,
such that the following property holds. Let $\Gamma_{s}^{t}$ be the unitary propagator generated by the linear Schr\"odinger equation \eqref{eq} in the time interval $[s,t]\subset[0,T]$ with the potential $V$ defined by \eqref{eq_potentiel_3}. For all $k\leq N$, there exists $\alpha_k\in\CC$ with $|\alpha_k|=1$ such that 
$$\big\|~\Gamma_{0}^{T}\sin(k\pi x)\,-\,\alpha_k\sin(\sigma(k)\pi x)~\big\|_{L^2}\leq \varepsilon~. $$
\end{prop}
\begin{demo}
Since the main arguments of the proof of Theorem \ref{th} still hold for $V$ given by \eqref{eq_potentiel_3}, it is enough to consider the ideal problem on the interval $[0,1]$ split in several subintervals. We show that any permutation $\sigma$ of $\NN^*$ can be realized by changing the locations $a_j(t)$ of the splittings. We examine the first $N$ eigenvalues. 
Let $$M=\max(N,\sigma(1),\ldots,\sigma(N)), \ \ \ \ \ \ J=M-1.$$ 

\vspace{3mm}

\noindent {\bf The initial ``vertical'' adiabatic motion.} At $t=0$, we split the interval $[0,1]$ in $M$ subintervals $[0,a_1(0)]\cup[a_1(0),a_2(0)]\cup \ldots \cup[a_J(0),1]$. We choose $a_j(0)$ such that the intervals have decreasing lengths which are very close. More precisely, we assume
$$|a_1(0)-0|~>~|a_2(0)-a_1(0)|~>~\ldots~>~|1-a_J(0)|~>~\frac12\, |a_1(0)-0|~.$$
When we increase adiabatically the intensity $I_j$ of the potential walls from $0$ to a large $I_*$, the first mode becomes localized in the first interval $[0,a_1(0)]$, the second mode in the second interval and so on. We notice that, since $|1-a_J(0)|>\frac12 |a_1(0)-0|$, the first mode corresponding to any subinterval corresponds to an eigenvalue lower than the one corresponding to the second mode of the first subinterval. 

\vspace{3mm}

\noindent {\bf The ``horizontal'' quasi-adiabatic motion.} We use the quasi-adiabatic motion introduced in the proof of Theorem \ref{th} (see Section \ref{section-proof-of-theorem}). Now, we change the location of the walls so that $[0,a_1(T)]$ has the $\sigma(1)-$th length among the intervals, $[a_1(T),a_2(T)]$ has the $\sigma(2)-$th length and so on until $[a_{N-1}(T),a_N(T)]$. For $j>N$, we proceed as follows. We make the length of $[a_{N}(T),a_{N+1}(T)]$ have order $p$ among the lengths of all the intervals, where $p$ is the least integer in $\NN^*\setminus\{\sigma(1),\ldots,\sigma(N)\}$, the interval 
$[a_{N+1}(T),a_{N+2}(T)]$ gets the $q-$th length with $q$ being the first integer in $\NN^*\setminus\{\sigma(1),\ldots,\sigma(N),p\}$, and so on. Moreover, we require (as above) that the ratio of any two lengths is strictly less than $2$. As a consequence, the second eigenvalue in the longest subinterval is higher than the first eigenvalue in the other subintervals. 

\vspace{3mm}

\noindent {\bf The final ``vertical'' adiabatic motion.} 
We adiabatically remove the walls by decreasing $I$ to $0$. The mode localized in the first subinterval corresponds to the $\sigma(1)-$th mode, i.e. $\sin(\sigma(1)\pi x)$, the mode localized in the second subinterval becomes $\sin(\sigma(2)\pi x)$, and so on. Thus, our quasi-adiabatic control has changed the modes $\sin(k\pi x)$ for $k\leq N$ into the modes $\sin(\sigma(k)\pi x)$, up to error terms as small as we want. 
\end{demo}

\subsection{Proof of Theorem \ref{th3}}\label{section_approx}
In this section, we use the construction of Proposition \ref{th2} in order to prove Theorem \ref{th3}.
It is enough to show that applying control of type (\ref{eq_potentiel_3}) we can move the initial state 
$A \sin(\pi \cdot)$ arbitrarily $L^2$-close to any given state $u$ with norm $A/\sqrt{2}$. By reversing time,
this would imply that starting in any $L^2$-neighborhood of $u$ we can drive the system to the state
$A \sin(\pi \cdot)$. Hence, since the Schr\"odinger equation preserves the $L^2$-norm, starting with $u$
we can drive the system to any $L^2$-neighborhood of $A \sin(\pi \cdot)$. Thus, given the initial and target states $u_i$ and $u_f$ of equal $L_2$-norm, we can first drive $u_i$ in the $\varepsilon/2$-neighborhood of
$A \sin(\pi \cdot)$ with $A=\sqrt{2} \|u_i\|_{L^2}$ and, after that, the result is driven to the 
$\varepsilon$-neighborhood of $u_f$ by exactly the same control which drives $A \sin(\pi \cdot)$ to the
$\varepsilon/2$-neighborhood of $u_f$.
\vspace{3mm}

{\noindent \bf Distributing the amplitudes.} So, let us show that given 
any set of real $c_k$, $k=1,\dots,N$, such that  $\sum c_k^2=1$ and any
$\alpha_k\in\UU=\{z\in\CC,|z|=1\}$, we can drive the initial state $\sin(\pi \cdot)$ arbitrarily close to
\begin{equation}\label{target}
u=\sum_{k=1}^N c_k\alpha_k\sin(k\pi \cdot)
\end{equation}
(since we work up to a small error term, we may choose $N$ large and assume that $u$ is only supported by the first $N$ modes).

\begin{figure}[ht]
\begin{center}
\resizebox{0.9\textwidth}{40mm}{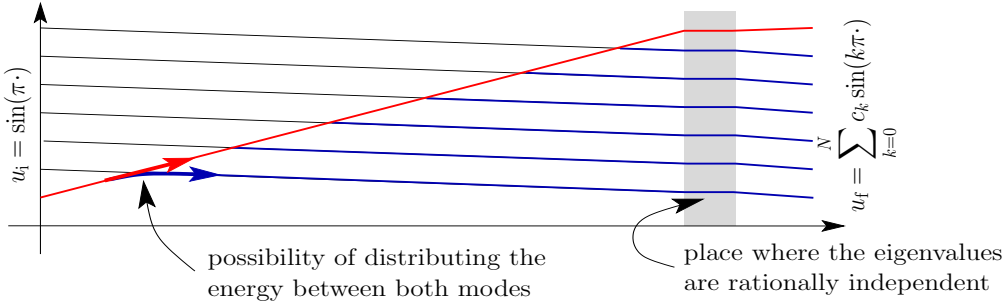}
\end{center}
\caption{\it The figure heuristically represents how the lines of eigenvalues behave when we distribute the amplitudes between the eigenmodes. In particular, if we start with $\ui=\sin(\pi\cdot)$, then it is possible to split the interval in $N$ subintervals and slowly modify their lengths to obtain the behavior illustrated above: the first eigenvalue crosses the next $N-1$ ones. At each crossing, we can choose the speed of the wall motion in the range between that described in Section \ref{section_fast_translation}, which would lead to complete crossing, and a much slower one which would make the system change adiabatically. This allows to create a superposition of the crossing modes with any given amplitudes. Moreover, if we make sure that the eigenvalues have irrational ratio at a particular time moment, then we can stop changing parameters for some time and
wait until the evolution of the autonomous Schr\"odinger equation sends the phases $\alpha_k$ to their target values.}\label{fig_th3}
\end{figure}

We start with the control described in Proposition \ref{th2} which enables to send the first mode to the $N-$th mode as in Figure \ref{fig_th3}, and modify the speed with which the potential walls move in order
to reach the vicinity of a state $u$ defined by (\ref{target}) with the given values of $c_k$ and arbitrary $\alpha_k$, $k=1,\dots, N$. We will tune the phases after that.

First, we split $[0,1]$ in $N$ subintervals with decreasing lengths, the largest subinterval being less
than twice longer than the smallest one. We adiabatically grow the potential walls at the boundary points between these intervals.
 
Second, we move each $a_j(t)$, $j=1, \dots, J$ to decrease the length of the first subinterval and to increase proportionally the length of the others. During this process, the line of eigenvalues starting from the lowest mode of the initial problem will cross the lines of eigenvalues corresponding to the modes 
$2,...,N$ of the initial problem, one by one, exactly in this order (see Figure \ref{fig_th3}; note that these
lines do not cross with each other). Before the crossing with the line of the second mode, the solution is close to the lowest energy mode $\psi_1$. By adiabatic theorem, this will remain the case if we move the walls extremely slowly (as the eigenvalues corresponding to $-\partial_{xx}^2+V$ with the potential \eqref{eq_potentiel_3} are simple, see Figure \ref{fig_crossing}). On the other hand, if we move the walls faster and make a short-time crossing as in Theorem \ref{th}, then after the crossing the solution will get close to the second (in the order of the increase of energy) mode $\psi_2$, see Section \ref{section_fast_translation}. By continuity, we can choose an intermediate speed in such a way that the solution after the crossing will be close to a superposition of these modes, namely to
$$c_1 \alpha_1 \psi_1 + \sqrt{1-c_1^2} \alpha_2 \psi_2$$
with some $\alpha_{1},\alpha_{2} \in \UU$ and $c_1$ given by (\ref{target}). Since the equation is linear, we can trace the evolution of the modes $\psi_1$ and $\psi_2$ separately. By construction, there will be no further
crossings for $\psi_1$, while here will be crossing of the line of eigenvalue corresponding to $\psi_2$ with the line corresponding to the third mode. As before, by tuning the speed of crossing, we can drive the mode $\psi_2$ to a superposition
$$\frac{c_2}{\sqrt{1-c_1^2}} \tilde\alpha_1 \psi_2 + \frac{\sqrt{1-c_1^2-c_2^2}}{\sqrt{1-c_1^2}} \tilde\alpha_2 \psi_3,$$
meaning that the original solution gets now close to
$$c_1 \alpha_1 \psi_1 + c_2 \alpha_2 \psi_2 + \sqrt{1-c_1^2-c_2^2} \alpha_3 \psi_3$$
with some new factors $\alpha_{1}$, $\alpha_2$ and $\alpha_{3}\in\UU$. One repeats this procedure until the last crossing
when $[0,a_1(T)]$ becomes the shortest interval. Then the solution gets closed to 
$$\sum_{k=1}^N c_k \alpha_k \psi_k$$
with exactly
the same coefficients $c_k$ as in (\ref{target}) and uncontrolled coefficients $\alpha_k\in\UU$.

\vspace{3mm}

{\noindent \bf Tuning the phases.} We can always arrange the above procedure in such a way that at the end of it
the lengths $L_1, \dots, L_N$ of the subintervals into which the potential walls divide the interval $(0,1)$ 
were rationally independent.
Since the walls are assumed to be very high, each of the first $N$ modes $\psi_1,\dots, \psi_N$ 
is close to the first eigenmode of the Laplacian with
the Dirichlet boundary conditions on some of the subintervals: the mode $\psi_k$ is mostly supported in 
the interval $[a_k,a_{k+1}]$ of length $L_k$ for $k=1,\dots, N-1$ and the mode $\psi_N$ is mostly supported in the interval
$[0,a_1]$ of the length $L_N$. The corresponding eigenvalues are close to $\lambda_k=\left(\frac{\pi}{L_k}\right)^2$, so
they are close to a set of rationally independent numbers.

Recall that we can make this approximation as good as we want if we choose $\eta_*$ and $I_*$ large enough (i.e., if we make
the potential wall sufficiently thin and high). Thus, if we just wait for a time $t$ in this situation, then the solution will get close to $$\sum_{k=1}^N c_k\alpha_k e^{-it\lambda_k} \psi_k.$$ 
Due to the rational independence, the flow $t\mapsto (e^{-it\lambda_k})_{k\leq N}$ is dense in the torus $\TT^N$. 
Thus, we can get the solution close to 
$$\sum_{k=1}^N c_k\hat\alpha_k \psi_k$$ 
with any $\hat\alpha_k\in\UU$ prescribed in advance. Note that the upper bound on the waiting time $t$ is independent of the 
initial and target values of $\alpha_k$, $k=1,\dots,N$, and depends only on the desired accuracy of the approximation.

Now we slowly decrease the height of the potential walls up to the extinction of the potential. This is an adiabatic process which
transforms (up to a small error) the modes $\psi_k$ into $\sin(\pi k \cdot)$ while keeping $c_k$ constant. The phases acquire
a shift, i.e., the solution gets close to 
$$\sum_{k=1}^N c_k\hat\alpha_ke^{i\theta_k} \psi_k$$ 
for some real $\theta_k$, $k=1,\dots,N$. Importantly, the phase shifts $\theta_k$ do not depend on $\hat\alpha_k$, so
choosing $\hat\alpha_k=e^{-i\theta_k}\alpha_k$ we obtain the desired result.



\end{document}